\documentclass[a4paper,10pt,reqno]{amsart}
\usepackage[a4paper,tmargin=40mm,bmargin=35mm,hmargin=30mm]{geometry}
\usepackage[T1]{fontenc}
\usepackage{lmodern}
\usepackage{amsmath}
\usepackage{amssymb}
\usepackage{amsthm}
\usepackage{stmaryrd}
\usepackage{tikz}
\usetikzlibrary{cd}
\usepackage{booktabs}
\usepackage{multirow}
\usepackage{multicol}
\usepackage{here}
\usepackage{aliascnt}
\usepackage{hyperref}
\usepackage[noabbrev]{cleveref}

\newcommand{\NewTheorem}[2]{
	\newaliascnt{#1}{TheoremEnvironment}
	\newtheorem{#1}[#1]{#1}
	\aliascntresetthe{#1}
	\crefname{#1}{#1}{#2}
	\Crefname{#1}{#1}{#2}
}

\theoremstyle{definition}

\NewTheorem{Definition}{Definitions}
\NewTheorem{Remark}{Remarks}
\NewTheorem{Axiom}{Axioms}
\NewTheorem{Example}{Examples}
\NewTheorem{Observation}{Observations}
\NewTheorem{Convention}{Conventions}
\NewTheorem{Notation}{Notations}
\NewTheorem{Setting}{Settings}
\NewTheorem{Question}{Questions}
\NewTheorem{Answer}{Answers}
\NewTheorem{Conjecture}{Conjectures}
\NewTheorem{Problem}{Problems}
\NewTheorem{Solution}{Solutions}
\NewTheorem{Goal}{Goals}
\NewTheorem{Comment}{Comments}
\NewTheorem{Aim}{Aims}
\NewTheorem{Caution}{Cautions}
\NewTheorem{Exercise}{Exercises}

\theoremstyle{plain}
\NewTheorem{Proposition}{Propositions}
\NewTheorem{Lemma}{Lemmas}
\NewTheorem{Theorem}{Theorems}
\NewTheorem{Corollary}{Corollaries}

\crefname{enumi}{}{}
\Crefname{enumi}{}{}
\creflabelformat{enumi}{(#2#1#3)}
\crefname{enumii}{}{}
\Crefname{enumii}{}{}
\creflabelformat{enumii}{(#2#1#3)}
\crefname{enumiii}{}{}
\Crefname{enumiii}{}{}
\creflabelformat{enumiii}{(#2#1#3)}

\makeatletter
\renewcommand{\p@enumii}{}
\renewcommand{\p@enumiii}{}
\makeatother

\numberwithin{equation}{section}
\crefname{equation}{}{}
\Crefname{equation}{}{}
\creflabelformat{equation}{(#2#1#3)}

\newcommand{\SwapSymbols}[1]{
	\expandafter\let\expandafter\temporarysymbol\csname #1\endcsname
	\expandafter\let\csname #1\expandafter\endcsname\csname var#1\endcsname
	\expandafter\let\csname var#1\endcsname\temporarysymbol
}

\SwapSymbols{epsilon}
\SwapSymbols{phi}
\SwapSymbols{Gamma}
\SwapSymbols{Delta}
\SwapSymbols{Theta}
\SwapSymbols{Lambda}
\SwapSymbols{Xi}
\SwapSymbols{Pi}
\SwapSymbols{Sigma}
\SwapSymbols{Upsilon}
\SwapSymbols{Phi}
\SwapSymbols{Psi}
\SwapSymbols{Omega}

\newcommand{\bbZ}{\mathbb{Z}}

\newcommand{\cA}{\mathcal{A}}

\newcommand{\cC}{\mathcal{C}}

\newcommand{\cE}{\mathcal{E}}

\newcommand{\cG}{\mathcal{G}}

\newcommand{\cM}{\mathcal{M}}
\newcommand{\cN}{\mathcal{N}}

\newcommand{\cS}{\mathcal{S}}

\newcommand{\cX}{\mathcal{X}}

\newcommand{\kp}{\mathfrak{p}}
\newcommand{\kq}{\mathfrak{q}}

\let\originalleft\left
\let\originalright\right
\renewcommand{\left}{\mathopen{}\mathclose\bgroup\originalleft}
\renewcommand{\right}{\aftergroup\egroup\originalright}

\newcommand{\into}{\hookrightarrow}

\newcommand{\onto}{\twoheadrightarrow}

\newcommand{\isoto}{\xrightarrow{\smash{\raisebox{-0.25em}{$\sim$}}}}
\newcommand{\isofrom}{\xleftarrow{\smash{\raisebox{-0.25em}{$\sim$}}}}

\newcommand{\set}[2][]{\mathopen{#1\{}#2\mathclose{#1\}}}

\newcommand{\setwithtext}[2][]{\mathopen{#1\{}\,\textnormal{#2}\,\mathclose{#1\}}}
\newcommand{\setwithcondition}[3][]{\mathopen{#1\{}\,#2\mathrel{#1|}#3\,\mathclose{#1\}}}

\newcommand{\op}{\textnormal{op}}

\DeclareMathOperator{\Hom}{Hom}
\DeclareMathOperator{\End}{End}

\DeclareMathOperator{\Ext}{Ext}

\DeclareMathOperator{\KGdim}{KGdim}
\DeclareMathOperator{\projdim}{proj.\!dim}
\DeclareMathOperator{\injdim}{inj.\!dim}
\DeclareMathOperator{\gldim}{gl.\!dim}

\DeclareMathOperator{\Set}{Set}

\DeclareMathOperator{\Pro}{Pro}
\DeclareMathOperator{\Ind}{Ind}

\DeclareMathOperator{\noeth}{noeth}

\DeclareMathOperator{\Mod}{Mod}
\let\mod\relax
\DeclareMathOperator{\mod}{mod}

\DeclareMathOperator{\Ann}{Ann}

\DeclareMathOperator{\Db}{\mathbf{D}^{\mathrm{b}}}

\DeclareMathOperator{\Zg}{Zg}

\DeclareMathOperator{\Spec}{Spec}

\DeclareMathOperator{\Ass}{Ass}
\DeclareMathOperator{\Supp}{Supp}

\DeclareMathOperator{\ASpec}{ASpec}

\DeclareMathOperator{\AAss}{AAss}
\DeclareMathOperator{\ASupp}{ASupp}

\DeclareMathOperator{\Epi}{Epi}
\DeclareMathOperator{\Const}{Const}
\DeclareMathOperator{\varExt}{\mathbb{E}xt}
\DeclareMathOperator{\cprojdim}{c.\!proj.\!dim}

\DeclareMathOperator{\DaHom}{D_{\alpha}Hom}
\DeclareMathOperator{\DaExt}{D_{\alpha}Ext}
\DeclareMathOperator{\varDaExt}{D_{\alpha}\mathbb{E}xt}

\newcommand{\resp}{resp.\ }

\title{Extension groups between atoms in abelian categories}

\author{Ryo Kanda}
\address{Department of Mathematics, Graduate School of Science, Osaka City University, 3-3-138, Sugimoto, Sumiyoshi, Osaka, 558-8585, Japan}
\email{ryo.kanda.math@gmail.com}

\begin{document}

\begin{abstract}
	We introduce the extension groups between atoms in an abelian category. For a locally noetherian Grothendieck category, the localizing subcategories closed under injective envelopes are characterized in terms of those extension groups. We also introduce the virtual duals of the extension groups between atoms to measure the global dimension of the category. A new topological property of atom spectra is revealed and it is used to relate the projective dimensions of atoms with the Krull-Gabriel dimensions. As a byproduct of the topological observation, we show that there exists a spectral space that is not homeomorphic to the atom spectrum of any abelian category.
\end{abstract}

\maketitle
\tableofcontents

\section{Introduction}
\label{30934094}

Classification of subcategories is one of the important problems and has been widely studied in several areas of mathematics. In the context of representation theory of rings, \emph{the} prototypical result was established by Gabriel:

\begin{Theorem}[{Gabriel \cite[Lemma~1 in p.~412 and Corollary~1 in p.~425]{MR0232821}}]\label{62164886}
	Let $R$ be a commutative noetherian ring. Then there is an order preserving bijection
	\begin{equation*}
		\setwithtext{localizing subcategories of $\Mod R$}\isoto\setwithtext{specialization-closed subsets of $\Spec R$}.
	\end{equation*}
\end{Theorem}
This result has been generalized in various ways; see \cite{MR1174255,MR1436741,MR2456904,MR2680200} and \cite[Proposition~4 in p.~446]{MR0232821}, for example. Among those generalizations, Herzog \cite{MR1434441} and Krause \cite{MR1426488} showed for a locally coherent Grothendieck category that there is an order-preserving bijection between the localizing subcategories of finite type and the open subsets of the Ziegler spectrum. The Ziegler spectrum is a topological space whose points are the isomorphism classes of indecomposable injective objects. Their results in particular imply the classification of all localizing subcategories for a locally noetherian Grothendieck category, which can be applied to the category of right modules over a right noetherian ring. On the other hand, the author observed that the Serre subcategories of a noetherian abelian category are classified in terms of the atom spectrum. The atom spectrum $\ASpec\cA$ of an abelian category $\cA$ is a topological space whose points, called \emph{atoms}, are the equivalence classes of monoform objects, and its topology is called the \emph{localizing topology}. The definition of atoms is based on work of Storrer \cite{MR0360717}. The atom spectrum is homeomorphic the Ziegler spectrum for a locally noetherian Grothendieck category, but it is still valid for a noetherian abelian category such as the category of finitely generated right modules over a right noetherian ring, and it allows us to show the aforementioned result.

Gabriel also showed a remarkable property of the category of modules over a commutative noetherian ring:

\begin{Theorem}[{Gabriel \cite[Proposition~10 in p.~428]{MR0232821}}]\label{66498719}
	Let $R$ be a commutative noetherian ring. Then every localizing subcategory of $\Mod R$ is closed under injective envelopes.
\end{Theorem}

For a locally noetherian Grothendieck category, localizing subcategories are not necessarily closed under injective envelopes (see \cref{23958094}). This means that \cref{66498719} uses some property that is specific to $\Mod R$ and it is natural to ask when localizing subcategories are closed under injective envelopes in general.

As mentioned above, the localizing subcategories of a locally noetherian Grothendieck category are in bijection with the (Ziegler or atom) spectrum. So one of the problems we should consider is: \emph{Characterize the localizing subcategories closed under injective envelopes in terms of the spectrum.} Papp considered this problem and gave several characterization (\cite{MR0393120,MR0409559,MR0506434}), but in this paper we take a different approach from those.

Our solution to this problem is given in terms of the \emph{extension groups between atoms}. Atoms in an abelian category $\cA$ can be regarded as pro-objects in $\cA$ (see \cref{15182006}) and we can define the extension groups $\Ext_{\cA}^{i}(\alpha,\beta)$ for atoms $\alpha,\beta\in\ASpec\cA$ in a natural way. One of our main results is the following:

\begin{Theorem}[\cref{07240397}]\label{92828160}
	Let $\cG$ be a locally noetherian Grothendieck category. Then there is an order-preserving bijective correspondence between
	\begin{itemize}
		\item localizing subcategories of $\cG$ that are closed under injective envelopes, and
		\item open subsets $\Phi$ of $\ASpec\cG$ with $\Ext_{\cG}^{1}(\alpha,\beta)=0$ for all $\alpha\in\ASpec\cG\setminus\Phi$ and $\beta\in\Phi$.
	\end{itemize}
\end{Theorem}

\cref{92828160} is one of the consequences of general observation for $\Ext_{\cG}^{i}(\alpha,\beta)$ for arbitrary integers $i\geq 0$. Since the extension groups between atoms are difficult to control due to its definition involving inverse limit, we also study a variant of them, which we call the \emph{virtual dual} of the extension groups and denote by $\varDaExt_{\cG}^{i}(\alpha,\beta)$. Indeed, we can determine the global dimension only using those virtual duals:

\begin{Theorem}[\cref{24385902}]\label{30418873}
	Let $\cG$ be a locally noetherian Grothendieck category. Then
	\begin{equation*}
		\gldim\cG=\sup\setwithcondition{i\geq 0}{\textnormal{$\varDaExt_{\cG}^{i}(\alpha,\beta)\neq 0$ for some $\alpha,\beta\in\ASpec\cG$}}.
	\end{equation*}
\end{Theorem}

The extension groups $\Ext_{\cG}^{i}(\alpha,N)$ between an atom $\alpha$ and an object $N$ have already been introduced in \cite{MR3272068}. There we gave a description of those for noetherian algebras. We will obtain a similar description for extension groups between atoms:

\begin{Theorem}[\cref{89835504}]\label{09558590}
	Let $\Lambda$ be a noetherian $R$-algebra. Let $i\geq 0$ be an integer and $P,Q\in\Spec\Lambda$. Then
	\begin{equation*}
		\Ext_{\Lambda}^{i}(\widetilde{P},\widetilde{Q})\cong
		\begin{cases}
			\Ext_{\Lambda_{\kp}}^{i}(S(P),S(Q)) & \text{if $P\cap R=Q\cap R=:\kp$,} \\
			0 & \text{otherwise,}
		\end{cases}
	\end{equation*}
	where $S(P)$ is the simple right $\Lambda_{\kp}$-module corresponding to $P$.
\end{Theorem}

The inverse systems that define the extension groups between atoms are often eventually constant, and in that case, we do not have to take the inverse limit. Indeed, in the setting of \cref{09558590}, the inverse systems are eventually constant when $P\cap R=Q\cap R$. We will seek such good cases for a locally noetherian Grothendieck category. We define the projective dimension $\projdim\alpha$ of an atom $\alpha$ in terms of vanishing of extension groups $\Ext_{\cG}^{i}(\alpha,-)$, and define its variant $\cprojdim\alpha$ to be the infimum of the integers $i$ such that the inverse limit defining $\Ext_{\cG}^{i}(\alpha,\beta)$ is eventually constant and nonzero. We will show that the difference of these two invariants of an atom is bounded by the Krull-Gabriel dimension of the category:

\begin{Theorem}[\cref{26816538}]\label{71731514}
	Let $\cG$ be a locally noetherian Grothendieck category. For every $\alpha\in\ASpec\cG$, we have
	\begin{equation*}
		\projdim\alpha\leq\cprojdim\alpha+\KGdim\cG.
	\end{equation*}
\end{Theorem}

The Krull-Gabriel dimension of a locally noetherian Grothendieck category is determined by the topology of the atom spectrum. In order to prove \cref{71731514}, we will show a new topological property of the atom spectrum of an abelian category. This also allows us to make an interesting observation.

Hochster \cite[Theorem~6]{MR0251026} showed that a topological space is homeomorphic to $\Spec R$, equipped with the Zariski topology, for a commutative ring $R$ if and only if the topological space is spectral (see \cref{60330213}). The topologies on atom spectra considered in this paper are not a generalization of Zariski topology. However, for a commutative noetherian ring $R$, the localizing topology on the atom spectrum of $\Mod R$ is the Hochster dual of the Zariski topology, which implies that the atom spectrum is also a spectral space. Although this is not necessarily true for a commutative ring in general, one would expect some connection between the topological spaces arising as atom spectra and spectral spaces. An abelian category is a massive generalization of the category of modules over a commutative ring, so a natural question would be: \emph{Is every spectral space homeomorphic to the atom spectrum of some abelian category?} Our topological observation implies that the answer is \emph{no}:

\begin{Theorem}[\cref{43496855}]\label{47194013}
	There exists a spectral space that is not homeomorphic to the atom spectrum of any abelian category equipped with the localizing topology.
\end{Theorem}

\subsection*{Acknowledgments}

The author would like to thank Hiroyuki Minamoto for stimulating discussion and the anonymous referee for their valuable comments. The author also thanks S.~Paul Smith for his hospitality as a host researcher at the University of Washington.

The author was a JSPS Overseas Research Fellow. This work was supported by JSPS KAKENHI Grant Numbers JP13J00249, JP16H06337, JP17K14164, and JP20K14288, Leading Initiative for Excellent Young Researchers, MEXT, Japan, and Osaka City University Advanced Mathematical Institute (MEXT Joint Usage/Research Center on Mathematics and Theoretical Physics JPMXP0619217849).

\section{Preliminaries}
\label{39845920}

\begin{Convention}\label{32985904}\leavevmode
	\begin{enumerate}
		\item\label{47350892} We fix a Grothendieck universe throughout the paper. A set is said to be \emph{small} if it belongs to the fixed universe. For a category, the collection of objects and that of morphisms are sets, which are not necessarily small. Every $\Hom$-set between two objects is supposed to be small. Every set arising as an index set of a colimit, a limit, or a generating set (see \cref{05592639} \cref{31409627}) should be in bijection with a small set. All rings and modules are assumed to be small.
		\item\label{23098044} Coproducts and products in an abelian category are called \emph{direct sums} and \emph{direct products}, respectively. A \emph{direct limit} means a colimit of a direct system indexed by a directed set. An \emph{inverse limit} means a limit of an inverse system indexed by a directed set. A directed (or inverse) system indexed by a directed set $I$ is often written as $\set{M_{i}}_{i\in I}$ by omitting the structure morphisms.
	\end{enumerate}
\end{Convention}

For a ring $\Lambda$, denote by $\Mod\Lambda$ the category of right $\Lambda$-modules. If $\Lambda$ is right noetherian, then denote by $\mod\Lambda$ the category of finitely generated right $\Lambda$-modules.

\subsection{Noetherian abelian categories}
\label{28374544}

First we recall the definitions of a generating set and a Grothendieck category:

\begin{Definition}\label{05592639}\leavevmode
	\begin{enumerate}
		\item\label{31409627} Let $\cA$ be an abelian category. A \emph{generating set} in $\cA$ is a set of objects $\set{U_{i}}_{i\in I}$ in $\cA$, where $I$ is in bijection with a small set, such that for every nonzero morphism $f\colon X\to Y$ in $\cA$, there exist $i\in I$ and a morphism $g\colon U_{i}\to X$ satisfying $fg\neq 0$.
		
		A \emph{cogenerating set} is a generating set in the opposite category.
		\item\label{50588886} A \emph{Grothendieck category} is an abelian category $\cG$ satisfying the following properties:
		\begin{itemize}
			\item $\cG$ admits direct sums (and hence all colimits).
			\item Direct limits are exact in $\cG$.
			\item $\cG$ admits a generating set.
		\end{itemize}
	\end{enumerate}
\end{Definition}

It is known that every Grothendieck category $\cG$ admits all limits and every object in $\cG$ has its injective envelope.

In most of the main results in this paper, we assume that a given abelian category satisfies a noetherian property. The noetherian property is one of the following two properties, depending on whether the abelian category is Grothendieck or skeletally small:

\begin{Definition}\label{63443637}\leavevmode
	\begin{enumerate}
		\item\label{41034765} A Grothendieck category $\cG$ is called \emph{locally noetherian} if it admits a generating set consisting of noetherian objects.
		\item\label{16823143} An abelian category $\cA$ is called \emph{noetherian} if all objects in $\cA$ are noetherian and $\cA$ is skeletally small, that is, the set of isomorphism classes in $\cA$ is in bijection with a small set.
	\end{enumerate}
\end{Definition}

\begin{Remark}\label{07365193}
	For a locally noetherian Grothendieck category $\cG$, an object $M$ in $\cG$ is noetherian if and only if $M$ is finitely generated (in the sense of \cite[section~V.3]{MR0389953}) if and only if $M$ is finitely presented (again in the sense of \cite[section~V.3]{MR0389953}). This is easily deduced from the definitions of those notions. Recall that $M$ is finitely presented if and only if the functor $\Hom_{\cG}(M,-)\colon\cG\to\Mod\bbZ$ preserves direct limits (\cite[Proposition~V.3.4]{MR0389953}).
\end{Remark}

The category $\Mod\Lambda$ for a right noetherian ring is an example of a locally noetherian Grothendieck category, and its full subcategory $\mod\Lambda$, which consists of all noetherian objects is a noetherian abelian category. The correspondence between these two categories is generalized as follows:

\begin{Theorem}[{\cite[Theorem~1 in p.~356]{MR0232821}}]\label{38744944}
	There is a bijective correspondence between
	\begin{itemize}
		\item equivalence classes of locally noetherian Grothendieck categories and
		\item equivalence classes of noetherian abelian categories.
	\end{itemize}
	Each locally noetherian Grothendieck category $\cG$ corresponds to its full subcategory
	\begin{equation*}
		\noeth\cG:=\setwithtext{noetherian objects in $\cG$}.
	\end{equation*}
\end{Theorem}

Now we want to show \cref{43825023}, which allows us to compute the global dimension of a locally noetherian Grothendieck category only using extension groups between noetherian objects. The proof uses the next two results.

\begin{Proposition}\label{06236705}
	Let $\cG$ be a locally noetherian Grothendieck category. Let $i\geq 0$ be an integer and $M,N\in\cG$. If $\Ext_{\cG}^{i}(M,N)\neq 0$, then there exists a noetherian subquotient $M'$ of $M$ such that $\Ext_{\cG}^{i}(M',N)\neq 0$.
\end{Proposition}

\begin{proof}
	Although the proof is similar to that for Baer's criterion (see \cite[2.3.1]{MR1269324}), we give a complete proof for the convenience of the reader.
	
	Assume $i=0$. Since $\cG$ is locally noetherian, the object $M$ is the direct limit of its noetherian subobjects: $M=\varinjlim_{i\in I}M_{i}$. Thus the claim follows from
	\begin{equation*}
		\Hom_{\cG}(\varinjlim_{i\in I}M_{i},N)\cong\varprojlim_{i\in I}\Hom_{\cG}(M_{i},N).
	\end{equation*}
	
	If $i\geq 2$, then we take a short exact sequence
	\begin{equation*}
		0\to N\to J\to N'\to 0,
	\end{equation*}
	where $J$ is an injective object. Then it induces an isomorphism
	\begin{equation*}
		\Ext_{\cG}^{i-1}(M,N')\isoto\Ext_{\cG}^{i}(M,N).
	\end{equation*}
	Repeating this, the problem is reduced to the case $i=1$.
	
	Let $i=1$. By $\Ext_{\cG}^{1}(M,N)\neq 0$, there exists a short exact sequence
	\begin{equation*}
		0\to N\to E\to M\to 0
	\end{equation*}
	that does not split. In other words, by regarding $N$ as a subobject of $E$, the identity morphism $N\to N$ cannot be extended to a morphism $E\to N$.
	
	Let $\cE$ be the set of pairs $(E',f)$, where $E'\subset E$ is a subobject containing $N$ and $f\colon E'\to N$ is a morphism whose restriction to $N$ is the identity. We define a partial order on $\cE$ by
	\begin{equation*}
		(E'_{1},f_{1})\leq (E'_{2},f_{2})\iff E'_{1}\subset E'_{2}\quad\text{and}\quad f_{2}|_{E'_{1}}=f_{1}.
	\end{equation*}
	For every totally ordered subset $\set{(E'_{i},f_{i})}_{i\in I}\subset\cE$, the direct limit
	\begin{equation*}
		\varinjlim_{i\in I}f_{i}\colon\varinjlim_{i\in I}E'_{i}\to N
	\end{equation*}
	gives its upper bound in $\cE$. Thus, by Zorn's lemma, $\cE$ has a maximal element $(E'_{0},f_{0})$.
	
	Since $f_{0}$ is an extension of the identity on $N$, we have $E'_{0}\subsetneq E$. Let $E'/E'_{0}$ be a nonzero noetherian subobject of $E/E'_{0}$. If $\Ext_{\cG}^{1}(E'/E'_{0},N)=0$, then the exact sequence
	\begin{equation*}
		\Hom_{\cG}(E',N)\to\Hom_{\cG}(E'_{0},N)\to\Ext_{\cG}^{1}(E'/E'_{0},N)=0
	\end{equation*}
	implies that $f_{0}\colon E'_{0}\to N$ can be extended to $E'\to N$. This contradicts the maximality of $(E'_{0},f_{0})$. Therefore $\Ext_{\cG}^{1}(E'/E'_{0},N)\neq 0$. Since $N\subset E'_{0}\subset E'\subset E$ and $E/N\cong M$, the object $M':=E'/E'_{0}$ is a subquotient of $M$.
\end{proof}

\begin{Proposition}\label{93588947}
	Let $\cG$ be a locally noetherian Grothendieck category. Let $i\geq 0$ be an integer and $M\in\cG$ a noetherian object. Then the functor
	\begin{equation*}
		\Ext_{\cG}^{i}(M,-)\colon\cG\to\Mod\bbZ
	\end{equation*}
	preserves direct limits.
\end{Proposition}

\begin{proof}
	This follows from \cite[Proposition~15.3.3]{MR2182076}. Indeed, $\Hom_{\cG}(M,-)$ is a left exact functor that preserves direct limits (see \cref{07365193}). The full subcategory of $\cG$ consisting of all injective objects is a $\Hom_{\cG}(M,-)$-injective additive subcategory by \cite[Corollary~13.3.8]{MR2182076}, and it is closed under direct limits since $\cG$ is locally noetherian (see \cite[Theorem~5.8.7]{MR0340375}). Thus \cite[Proposition~15.3.3]{MR2182076} is applicable.
\end{proof}

For an abelian category $\cA$, its \emph{global dimension} is
\begin{equation*}
	\gldim\cA:=\sup\setwithcondition{i\geq 0}{\textnormal{$\Ext^{i}_{\cA}(M,N)\neq 0$ for some $M,N\in\cA$}}.
\end{equation*}
Let $\gldim\cA:=-1$ if $\cA$ is zero.

The \emph{injective dimension} of an object $N\in\cA$ is defined to be
\begin{equation*}
	\injdim N:=\sup\setwithcondition{i\geq 0}{\textnormal{$\Ext^{i}_{\cA}(M,N)\neq 0$ for some $M\in\cA$}}.
\end{equation*}
Let $\injdim N:=-1$ if $N=0$.

\begin{Proposition}\label{43825023}
	Let $\cG$ be a locally noetherian Grothendieck category. Then
	\begin{equation*}
		\gldim\cG=\gldim(\noeth\cG).
	\end{equation*}
\end{Proposition}

\begin{proof}
	This follows from \cref{06236705,93588947}.
\end{proof}

\subsection{Atom spectrum}
\label{30980844}

The atom spectrum of an abelian category is the main object to study in this paper. We recall its definition and some basic properties. For further results on atom spectra, see \cite{arXiv:1711.06946}, for example.

\begin{Definition}\label{82443509}
	Let $\cA$ be an abelian category.
	\begin{enumerate}
		\item\label{34549823} A \emph{monoform object} in $\cA$ is a nonzero object $H\in\cA$ that has no subobjects $0\neq L'\subsetneq L\subset H$ and $0\neq N\subset H$ such that $L/L'\cong N$.
		\item\label{34958904} We say that two monoform objects $H_{1}$ and $H_{2}$ are \emph{atom-equivalent} if there exist nonzero subobjects $L_{1}\subset H_{1}$ and $L_{2}\subset H_{2}$ such that $L_{1}\cong L_{2}$.
		\item\label{34989049} The \emph{atom spectrum} of $\cA$ is defined to be
		\begin{equation*}
			\ASpec\cA:=\frac{\setwithtext{monoform objects in $\cA$}}{\text{atom-equivalence}}.
		\end{equation*}
		An element of $\ASpec\cA$ is called an \emph{atom} in $\cA$. For each monoform object $H$, its equivalence class is denoted by $\overline{H}$.
	\end{enumerate}
\end{Definition}

\begin{Remark}\label{82094844}
	Let $\cA$ be an abelian category.
	\begin{enumerate}
		\item\label{09804444} A nonzero object $U\in\cA$ is called \emph{uniform} if any nonzero subobjects $L_{1}$ and $L_{2}$ of $U$ have nonzero intersection in $U$. Every monoform object is uniform (\cite[Proposition~2.6]{MR2964615}). This implies that the atom-equivalence is an equivalence relation between monoform objects (\cite[Proposition~2.8]{MR2964615}).
		\item\label{92837400} Every nonzero subobject of a monoform (\resp uniform) object is again monoform (\resp uniform) (\cite[Proposition~2.2]{MR2964615}).
		\item\label{09890449} The atom spectrum of $\cA$ is in bijection with a small set if $\cA$ admits a generating set (which is indexed by a small set; the proof of \cite[Proposition~2.7 (2)]{MR3351569} using \cite[Proposition~IV.6.6]{MR0389953} works also for any abelian category with a (small) generating set). Later we will focus on Grothendieck categories, for which these properties are satisfied.
		\item\label{17457124} Every nonzero noetherian object in $\cA$ has a monoform subobject (\cite[Theorem~2.9]{MR2964615}). Hence, if $\cA$ admits a generating set consisting of noetherian objects, then every nonzero object in $\cA$ has a monoform subobject.
	\end{enumerate}
\end{Remark}

The atom spectrum of an abelian category can be regarded as a generalization of
\begin{itemize}
	\item the set of prime ideals of a commutative ring (\cref{30989084}),
	\item the underlying space of a locally noetherian scheme (\cite[Theorem~7.6]{MR3452186}),
	\item the set of isomorphism classes of simple right modules over a right artinian ring (\cite[Proposition~8.2]{MR2964615}), and
	\item the set of prime two-sided ideals of a noetherian algebra (\cref{63482658}).
\end{itemize}
The next definition gives a generalized notion of associated points and supports:

\begin{Definition}\label{94890423}
	Let $\cA$ be an abelian category and let $M\in\cA$ be an object.
	\begin{enumerate}
		\item\label{98793433} The set of \emph{associated atoms} of $M$ is defined to be
		\begin{equation*}
			\AAss M:=\setwithcondition{\overline{H}\in\ASpec\cA}{\textnormal{$H$ is a monoform subobject of $M$}}.
		\end{equation*}
		\item\label{89840944} The \emph{atom support} of $M$ is defined to be
		\begin{equation*}
			\ASupp M:=\setwithcondition{\overline{H}\in\ASpec\cA}{\textnormal{$H$ is a monoform subquotient of $M$}}.
		\end{equation*}
	\end{enumerate}
\end{Definition}

\begin{Remark}\label{98409844}
	Associated atoms and atom supports are compatible with short exact sequences, direct unions, and direct sums, in the way that associated primes and supports of modules over commutative rings are (\cite[Proposition~2.6]{arXiv:1711.06946}). For example, if $0\to L\to M\to N\to 0$ is a short exact sequence in an abelian category $\cA$, then we have
	\begin{equation*}
		\AAss L\subset\AAss M\subset\AAss L\cup\AAss N
	\end{equation*}
	and
	\begin{equation*}
		\ASupp M=\ASupp L\cup\ASupp N.
	\end{equation*}
\end{Remark}

\begin{Definition}\label{30598443}
	Let $\cG$ be a Grothendieck category and let $\alpha=\overline{H}\in\ASpec\cG$. Define the isomorphism class of the \emph{injective envelope} $E(\alpha)$ of $\alpha$ to be the isomorphism class of the injective envelope $E(H)$ of $H$.
\end{Definition}

\begin{Remark}\label{30980944}
	We do not specify the representative $H$ in \cref{30598443} because the isomorphism class of $E(\alpha)$ does not depend on the choice of $H$. Indeed, for two monoform objects $H$ and $H'$ with $\alpha=\overline{H}=\overline{H'}$, their injective envelopes $E(H)$ and $E(H')$ are isomorphic to each other (\cite[Lemma~5.8]{MR2964615}), but there is no canonical isomorphism. We have a monomorphism $H\into E(H)\cong E(\alpha)$, but we do not have a canonical embedding $H\into E(\alpha)$.
\end{Remark}

\begin{Remark}\label{30989844}
	If $\cG$ is a locally noetherian Grothendieck category, then the correspondence $\alpha\mapsto E(\alpha)$ gives a bijection between $\ASpec\cG$ and the set $\Zg\cG$ of isomorphism classes of indecomposable injective objects in $\cG$ (\cite[Theorem~5.9]{MR2964615}). The set $\Zg\cG$ together with a certain topology is called the \emph{Ziegler spectrum} of $\cG$ (see \cite[Theorem~3.4]{MR1434441} or \cite[Lemma~4.1]{MR1426488} for a locally coherent Grothendieck category, \cite[Definition~5.7]{MR2964615} for the special case of a locally noetherian Grothendieck category). The bijection is in fact a homeomorphism for a locally noetherian Grothendieck category (\cite[Theorem~5.9]{MR2964615}).
	
	When we consider extension groups between atoms, it is more suitable to use atoms rather than indecomposable injective objects. Indeed, atoms will be regarded as pro-objects in \cref{15182006}, and if an atom is represented by a simple object, then the atom is isomorphic to the simple object as a pro-object (\cref{10271676}), but not to the corresponding indecomposable injective object, in general. Moreover, we can state our results also for a noetherian abelian category in terms of atoms, without mentioning the corresponding locally noetherian Grothendieck category.
\end{Remark}

We always assume that the following topology and partial order are defined on atom spectra:

\begin{Definition}\label{40598094}
	Let $\cA$ be an abelian category.
	\begin{enumerate}
		\item\label{20948044} There is a topology on $\ASpec\cA$ such that
		\begin{equation*}
			\setwithcondition{\ASupp M}{M\in\cA}
		\end{equation*}
		is an open basis (\cite[Proposition~3.2]{MR3351569}). We call it the \emph{localizing topology}.
		\item\label{32989044} Define a binary relation $\leq$ on $\ASpec\cA$, which is called the \emph{specialization order}, by
		\begin{equation*}
			\alpha\leq\beta\quad\text{if and only if}\quad\alpha\in\overline{\set{\beta}},
		\end{equation*}
		where $\overline{\set{\beta}}$ is the closure of the singleton $\set{\beta}$ with respect to the localizing topology.
	\end{enumerate}
\end{Definition}

\begin{Remark}\label{43980984}
	Let $\cA$ be an abelian category.
	\begin{enumerate}
		\item\label{60221442} A subset $\Phi\subset\ASpec\cA$ is open if and only if every $\alpha\in\Phi$ admits a monoform object $H\in\cA$ such that $\overline{H}=\alpha$ and $\ASupp H\subset\Phi$ (\cite[Definition~3.1 and Proposition~3.2]{MR3351569}).
		\item\label{29387433} The atom spectrum $\ASpec\cA$ is known to be a \emph{Kolmogorov space} (also called a \emph{$T_{0}$-space}; \cite[Proposition~3.5]{MR3351569}), that is, for any two distinct points $\alpha\neq\beta$ in $\ASpec\cA$, there exists an open subset $\Phi\subset\ASpec\cA$ such that $\set{\alpha,\beta}\cap\Phi$ consists of exactly one point. In other words, any two points are topologically distinguishable.
		\item\label{30984488} For a topological space, the binary relation defined in the way of \cref{40598094} \cref{32989044} is called the \emph{specialization preorder}, which is in general a partial preorder. Since $\ASpec\cA$ is a Kolmogorov space, the relation $\leq$ is a partial order.
	\end{enumerate}
\end{Remark}

Contrary to the Ziegler spectrum, the atom spectrum is defined even for a noetherian abelian category, and it is naturally identified with that of the corresponding locally noetherian Grothendieck category in the sense of \cref{38744944}:

\begin{Proposition}\label{89431240}
	Let $\cG$ be a locally noetherian Grothendieck category. Then there is a homeomorphism
	\begin{equation*}
		\ASpec(\noeth\cG)\isoto\ASpec\cG
	\end{equation*}
	given by $\overline{H}\mapsto\overline{H}$.
\end{Proposition}

\begin{proof}
	\cite[Proposition~5.3]{MR2964615}.
\end{proof}

\begin{Remark}\label{30989084}
	The atom spectrum of an abelian category is a generalization of the prime spectrum of a commutative ring. More precisely, the following assertions hold for every commutative ring $R$:
	\begin{enumerate}
		\item\label{20980444} (\cite[p.~631]{MR0360717}) There is a bijection
		\begin{equation*}
			\Spec R\isoto\ASpec(\Mod R)
		\end{equation*}
		given by $\kp\mapsto\overline{R/\kp}$.
		\item\label{03989044} (\cite[Proposition~2.13]{MR3351569}) For every $R$-module $M$, the bijection in \cref{20980444} induces bijections
		\begin{equation*}
			\Ass_{R}M\isoto\AAss M\quad\text{and}\quad\Supp_{R}M\isoto\ASupp M.
		\end{equation*}
		\item\label{30980948} (\cite[Proposition~7.2 (2)]{MR2964615}) A subset $\Phi\subset\ASpec(\Mod R)$ is open with respect to the localizing topology if and only if the inverse image of $\Phi$ by the bijection in \cref{20980444} is specialization-closed, that is, whenever it contains a prime ideal $\kp\subset R$, it also contains all prime ideals larger than $\kp$.
		\item\label{09809288} (\cite[Proposition~4.3]{MR3351569}) The bijection in \cref{20980444} is an isomorphism of partially ordered sets:
		\begin{equation*}
			(\Spec R,\subset)\isoto(\ASpec(\Mod R),\leq).
		\end{equation*}
	\end{enumerate}
\end{Remark}

\subsection{Serre subcategories and localizing subcategories}
\label{40980944}

We recall the definitions of a Serre subcategory and a localizing subcategory and state some fundamental results including the relation to atom spectra.

\begin{Definition}\label{20949889}\leavevmode
	\begin{enumerate}
		\item\label{30958094} Let $\cA$ be an abelian category. A full subcategory $\cS$ of $\cA$ is called a \emph{Serre subcategory} if it is closed under subobjects, quotient objects, and extensions, or equivalently: for every short exact sequence
		\begin{equation*}
			0\to L\to M\to N\to 0
		\end{equation*}
		in $\cA$, the object $M$ belongs to $\cS$ if and only if both $L$ and $N$ belong to $\cS$.
		\item\label{10938484} Let $\cG$ be a Grothendieck category. A Serre subcategory $\cX$ of $\cG$ is called a \emph{localizing subcategory} if it is moreover closed under direct sums.
	\end{enumerate}
\end{Definition}

\begin{Remark}\label{28742849}
	We summarize some facts on quotient categories here. See \cite[Chapter~III]{MR0232821} or \cite[section~4]{MR0340375} for more details.
	
	If $\cS$ is a Serre subcategory of an abelian category $\cA$, we can form the abelian category $\cA/\cS$ called the \emph{quotient category} together with a canonical functor $F\colon\cA\to\cA/\cS$, which is dense and exact. An object in $\cA$ is sent to zero by $F$ if and only if it belongs to $\cS$. For every object $M\in\cA$ and every subobject $L'\subset F(M)$, there exists a subobject $L\subset M$ such that $F(L)=L'$ as a subobject of $F(M)$ (\cite[Corollary~1 in p.~368]{MR0232821}).
	
	$\cS$ is called a \emph{localizing subcategory} if the canonical functor $\cA\to\cA/\cS$ admits a right adjoint. This definition agrees with \cref{20949889} \cref{10938484} when $\cA$ is a Grothendieck category (\cite[Proposition~4.6.3]{MR0340375}). If $\cA$ is a Grothendieck category and $\cS\subset\cA$ is a localizing subcategory, then $\cA/\cS$ is again a Grothendieck category (\cite[Proposition~9 in p.~378]{MR0232821}).
\end{Remark}

\begin{Proposition}\label{25716654}
	Let $\cG$ be a locally noetherian Grothendieck category.
	\begin{enumerate}
		\item\label{45393500} There is an order-preserving bijection
		\begin{equation*}
			\setwithtext{localizing subcategories of $\cG$}\isoto\setwithtext{Serre subcategories of $\noeth\cG$}
		\end{equation*}
		given by $\cX\mapsto\cX\cap\noeth\cG$. The inverse map sends each Serre subcategory $\cS\subset\noeth\cG$ to the smallest localizing subcategory of $\cG$ containing $\cS$.
		\item\label{05256887} Let $\cX$ be a localizing subcategory of $\cG$. Then $\cG/\cX$ is a locally noetherian Grothendieck category, and the inclusion $\noeth\cG\into\cG$ induces an equivalence
		\begin{equation*}
			\frac{\noeth\cG}{\cX\cap\noeth\cG}\isoto\noeth\frac{\cG}{\cX}.
		\end{equation*}
	\end{enumerate}
\end{Proposition}

\begin{proof}
	\cref{45393500} \cite[Proposition~10 in p.~379]{MR0232821}.
	
	\cref{05256887} This is a special case of \cite[Theorem~2.6]{MR1426488}.
\end{proof}

The following operations relate subcategories of a given abelian category and subsets of its atom spectrum:

\begin{Definition}\label{04958908}
	Let $\cA$ be an abelian category.
	\begin{enumerate}
		\item\label{92338444} For a Serre subcategory $\cS\subset\cA$, define the open subset $\ASupp\cS\subset\ASpec\cA$ by
		\begin{equation*}
			\ASupp\cS:=\bigcup_{M\in\cS}\ASupp M.
		\end{equation*}
		\item\label{83409202} For an open subset $\Phi\subset\ASpec\cA$, define the Serre subcategory $\ASupp^{-1}\Phi\subset\cA$ by
		\begin{equation*}
			\ASupp^{-1}\Phi:=\setwithcondition{M\in\cA}{\ASupp M\subset\Phi}.
		\end{equation*}
	\end{enumerate}
\end{Definition}

\begin{Theorem}[{Herzog \cite[Theorem~3.8]{MR1434441}, Krause \cite[Theorem~4.2]{MR1426488}, and Kanda \cite[Theorem~5.5]{MR2964615}}]\label{09809844}\leavevmode
	\begin{enumerate}
		\item\label{93450984} Let $\cG$ be a locally noetherian Grothendieck category. Then there is an order-preserving bijection
		\begin{equation*}
			\setwithtext{localizing subcategories of $\cG$}\isoto\setwithtext{open subsets of $\ASpec\cG$}
		\end{equation*}
		given by $\cX\mapsto\ASupp\cX$ whose inverse map is $\Phi\mapsto\ASupp^{-1}\Phi$.
		\item\label{42375944} Let $\cA$ be a noetherian abelian category. Then there is an order-preserving bijection
		\begin{equation*}
			\setwithtext{Serre subcategories of $\cA$}\isoto\setwithtext{open subsets of $\ASpec\cA$}
		\end{equation*}
		given by $\cS\mapsto\ASupp\cS$ whose inverse map $\Phi\mapsto\ASupp^{-1}\Phi$.
	\end{enumerate}
\end{Theorem}

\begin{Remark}\label{55074566}
	If $\cG$ is a locally noetherian Grothendieck category and $\cA=\noeth\cG$, then we have the commutative diagram
	\begin{equation*}
		\begin{tikzcd}
			\setwithtext{localizing subcategories of $\cG$}\ar[d,"\wr"']\ar[r,"\sim"] & \setwithtext{Serre subcategories of $\cA$}\ar[d,"\wr"] \\
			\setwithtext{open subsets of $\ASpec\cG$}\ar[r,"\sim"'] & \setwithtext{open subsets of $\ASpec\cA$}\rlap{,}
		\end{tikzcd}
	\end{equation*}
	where the bottom horizontal bijection is induced from the homeomorphism $\ASpec\cA\isoto\ASpec\cG$ in \cref{89431240}.
\end{Remark}

Atom spectra are compatible with taking quotient categories:

\begin{Theorem}\label{09609789}\leavevmode
	\begin{enumerate}
		\item\label{38104494} Let $\cG$ be a Grothendieck category and let $\cX\subset\cG$ be a localizing subcategory. Then there is a homeomorphism
		\begin{equation*}
			\ASpec\cG\setminus\ASupp\cX\isoto\ASpec\frac{\cG}{\cX}
		\end{equation*}
		given by $\overline{H}\mapsto\overline{F(H)}$, where $F\colon\cG\to\cG/\cX$ is the canonical functor.
		\item\label{52401192} Let $\cA$ be a noetherian abelian category and let $\cS\subset\cA$ be a Serre subcategory. Then there is a homeomorphism
		\begin{equation*}
			\ASpec\cA\setminus\ASupp\cS\isoto\ASpec\frac{\cA}{\cS}
		\end{equation*}
		given by $\overline{H}\mapsto\overline{F(H)}$, where $F\colon\cA\to\cA/\cS$ is the canonical functor.
	\end{enumerate}
\end{Theorem}

\begin{proof}
	\cref{38104494} \cite[Theorem~5.17]{MR3351569}.
	
	\cref{52401192} This is a combination of \cref{38104494}, \cref{55074566}, and \cref{25716654} \cref{05256887}.
\end{proof}

\begin{Remark}\label{26229957}
	The noetherian assumption in \cref{09609789} \cref{52401192} cannot be dropped.
	
	To see this, consider the category $\Mod^{\bbZ}k[x]$ of $\bbZ$-graded $k[x]$-modules whose morphisms are degree-preserving homomorphisms, where $k$ is a field and $k[x]$ is graded as $\deg x=1$. For a graded module $M=\bigoplus_{i\in\bbZ}M_{i}\in\Mod^{\bbZ}k[x]$, its degree shifts are denoted by $M(j)$ ($j\in\bbZ$) with the convention $M(j)_{i}=M_{i+j}$. Denote by $\mod^{\bbZ}k[x]$ its full subcategory consisting of all finitely generated modules. The category $\Mod^{\bbZ}k[x]$ is a locally noetherian Grothendieck category since the set $\setwithcondition{k[x](j)}{j\in\bbZ}$ is a generating set consisting of noetherian objects, and $\noeth(\Mod^{\bbZ}k[x])=\mod^{\bbZ}k[x]$ by \cref{07365193} (finite generation in $\Mod^{\bbZ}k[x]$ is the usual finite generation of modules). So \cref{89431240} implies that
	\begin{equation*}
		\ASpec(\Mod^{\bbZ}k[x])\setminus\ASupp(\mod^{\bbZ}k[x])=\emptyset.
	\end{equation*}
	Recall that the \emph{Matlis dual} of a $\bbZ$-graded $k[x]$-module $M=\bigoplus_{i\in\bbZ}M_{i}$ is defined to be
	\begin{equation*}
		M^{\vee}:=\bigoplus_{i\in\bbZ}\Hom_{k}(M_{-i},k)\in\Mod^{\bbZ}k[x].
	\end{equation*}
	Let $I:=k[x]^{\vee}\in\Mod^{\bbZ}k[x]$. It is easy to check that all proper subobjects of $I$ are of the form $(k[x]/x^{n}k[x])^{\vee}$ for $n\in\bbZ_{\geq 0}$, which all belong to $\mod^{\bbZ}k[x]$. Thus, by \cref{28742849}, the module $I$ is sent to a simple object by the canonical functor $\Mod^{\bbZ}k[x]\to\Mod^{\bbZ}k[x]/\mod^{\bbZ}k[x]$. Since the simple object is monoform, we obtain
	\begin{equation*}
		\ASpec\frac{\Mod^{\bbZ}k[x]}{\mod^{\bbZ}k[x]}\neq\emptyset.
	\end{equation*}
	Thus \cref{09609789} \cref{52401192} does not apply to the abelian category $\cA:=\Mod^{\bbZ}k[x]$ and its Serre subcategory $\cS:=\mod^{\bbZ}k[x]$.
\end{Remark}

\subsection{Krull-Gabriel dimension}
\label{26602600}

We recall the definition of the Krull-Gabriel dimensions of Grothendieck categories and their objects, which generalizes the Krull dimension of a commutative noetherian ring (see \cref{62636912}). We also define the Krull-Gabriel dimensions of atoms in a natural way.

\begin{Definition}[{\cite[Chapter~IV.1]{MR0232821}}]\label{13578089}
	Let $\cG$ be a Grothendieck category.
	\begin{enumerate}
		\item\label{14422843} For ordinal numbers $\lambda$ and $\lambda=-1$, we define the localizing subcategories $\cG_{\lambda}\subset\cG$ inductively as follows:
		\begin{itemize}
			\item $\cG_{-1}$ consists of all zero objects in $\cG$.
			\item $\cG_{\lambda+1}$ is the smallest localizing subcategory of $\cG$ containing all objects $M\in\cG$ that are sent to objects of finite length by the canonical functor $\cG\to\cG/\cG_{\lambda}$.
			\item For a limit ordinal $\lambda$, $\cG_{\lambda}$ is the smallest localizing subcategory of $\cG$ containing $\cG_{\mu}$ for all $\mu<\lambda$.
		\end{itemize}
		\item\label{26228506} For an object $M\in\cG$, its \emph{Krull-Gabriel dimension} is defined to be
		\begin{equation*}
			\KGdim M:=\inf\setwithcondition{\lambda}{M\in\cG_{\lambda}}.
		\end{equation*}
		\item\label{52765482} The \emph{Krull-Gabriel dimension} of $\cG$ is defined to be
		\begin{equation*}
			\KGdim\cG:=\inf\setwithcondition{\lambda}{\cG_{\lambda}=\cG}.
		\end{equation*}
	\end{enumerate}
	If the set in the definition of \cref{26228506} (\resp \cref{52765482}) is empty, then we say that \emph{the Krull-Gabriel dimension of $M$ (\resp $\cG$) does not exist}.
\end{Definition}

\begin{Remark}\label{80274336}
	It is known that the Krull-Gabriel dimension exists for every locally noetherian Grothendieck category (\cite[Proposition~7 in p.~387]{MR0232821}).
\end{Remark}

\begin{Remark}\label{69279173}
	We can define the Krull-Gabriel dimension of a noetherian abelian category $\cA$ analogously. For ordinal numbers $\lambda$ and $\lambda=-1$, define the Serre subcategories $\cA_{\lambda}\subset\cA$ inductively as follows:
	\begin{itemize}
		\item $\cA_{-1}$ consists of all zero objects in $\cA$.
		\item $\cA_{\lambda+1}$ is the Serre subcategory of $\cA$ consisting of all objects $M\in\cA$ that are sent to objects of finite length by the canonical functor $\cA\to\cA/\cA_{\lambda}$.
		\item For a limit ordinal $\lambda$, $\cA_{\lambda}$ is the union of all $\cA_{\mu}$ with $\mu<\lambda$.
	\end{itemize}
	We define the \emph{Krull-Gabriel dimension} of an object $M\in\cA$ to be
	\begin{equation*}
		\KGdim M:=\inf\setwithcondition{\lambda}{M\in\cA_{\lambda}}
	\end{equation*}
	and define the \emph{Krull-Gabriel dimension} of $\cA$ to be
	\begin{equation*}
		\KGdim\cA:=\inf\setwithcondition{\lambda}{\cA_{\lambda}=\cA}.
	\end{equation*}
	
	If $\cG$ is a locally noetherian Grothendieck category satisfying $\noeth\cG=\cA$, then \cref{25716654} implies that $\cA_{\lambda}=\cG_{\lambda}\cap\cA$ for all ordinal numbers $\lambda$ and $\lambda=-1$. Thus, for every object $M\in\cA$, its Krull-Gabriel dimension defined in $\cA$ is equal to that defined in $\cG$. We also have
	\begin{equation*}
		\KGdim\cA=\KGdim\cG.
	\end{equation*}
\end{Remark}

For a locally noetherian Grothendieck category, there is an order-preserving bijective correspondence between the localizing subcategories and the open subsets of the atom spectrum (\cref{09809844}). We will see that the Krull-Gabriel dimension can also be defined using the topological structure of the atom spectrum (\cref{63990180}).

\begin{Definition}\label{38948897}
	For a topological space $X$, define the open subspaces $X_{\lambda}$ for ordinal numbers $\lambda$ and $\lambda=-1$ inductively as follows:
	\begin{itemize}
		\item $X_{-1}=\emptyset$.
		\item $X_{\lambda+1}=X_{\lambda}\cup\setwithtext{open points of $X\setminus X_{\lambda}$}$, where an open point of $X\setminus X_{\lambda}$ means a point $x\in X\setminus X_{\lambda}$ such that $\set{x}$ is an open subset of the topological space $X\setminus X_{\lambda}$.
		\item For a limit ordinal $\lambda$, $X_{\lambda}$ is the union of all $X_{\mu}$ with $\mu<\lambda$.
	\end{itemize}
\end{Definition}

\begin{Definition}\label{13750276}
	Let $\cA$ be a locally noetherian Grothendieck category or a noetherian abelian category. The \emph{Krull-Gabriel dimension} of $\alpha\in\ASpec\cA$ is defined to be
	\begin{equation*}
		\KGdim\alpha:=\inf\setwithcondition{\lambda}{\alpha\in(\ASpec\cA)_{\lambda}}.
	\end{equation*}
\end{Definition}

\begin{Proposition}\label{60891414}
	Let $\cA$ be a locally noetherian Grothendieck category or a noetherian abelian category. For all ordinal numbers $\lambda$ and $\lambda=-1$, we have
	\begin{equation*}
		\ASupp(\cA_{\lambda})=\setwithcondition{\alpha\in\ASpec\cA}{\KGdim\alpha\leq\lambda}.
	\end{equation*}
\end{Proposition}

\begin{proof}
	Let $\cA$ be a noetherian abelian category. We use induction on $\lambda$. Note that the right-hand side of the equation is $(\ASpec\cA)_{\lambda}$.
	
	If $\lambda=-1$, then the both sides of the equation are empty.
	
	Let $\lambda$ be an arbitrary ordinal number. By the induction hypothesis, the desired equation for $\lambda+1$ follows once we prove
	\begin{equation*}
		\ASupp(\cA_{\lambda+1})\setminus\ASupp(\cA_{\lambda})=\setwithtext{open points of $\ASpec\cA\setminus\ASupp(\cA_{\lambda})$}.
	\end{equation*}
	Denote by $F\colon\cA\to\cA/\cA_{\lambda}$ the canonical functor.
	
	Let $\alpha\in\ASupp(\cA_{\lambda+1})\setminus\ASupp(\cA_{\lambda})$. Then $\alpha=\overline{H}$ for some monoform object $H\in\cA_{\lambda+1}$. Since $\alpha\notin\ASupp(\cA_{\lambda})$, \cref{09609789} implies that $F(H)$ is a monoform object. Since $F(H)$ is of finite length, the atom $\overline{F(H)}\in\ASpec(\cA/\cA_{\lambda})$ is represented by a simple subobject of $F(H)$. By \cite[Proposition~3.7 (1)]{MR3351569}, $\overline{F(H)}$ is an open point of $\ASpec(\cA/\cA_{\lambda})$. Again by \cref{09609789}, $\alpha=\overline{H}$ is an open point of $\ASpec\cA\setminus\ASupp(\cA_{\lambda})$.
	
	Conversely, let $\alpha=\overline{H}\in\ASpec\cA\setminus\ASupp(\cA_{\lambda})$ be an open point. Since $\overline{F(H)}\in\ASpec(\cA/\cA_{\lambda})$ is an open point, it is represented by a simple object $S$, again by \cite[Proposition~3.7 (1)]{MR3351569}. Since $F(H)$ and $S$ are atom-equivalent, we can regard $S$ as a subobject of $F(H)$. As in \cref{28742849}, there exists a subobject $H'\subset H$ such that $F(H')=S$ as a subobject of $F(H)$. Since $H'\in\cA_{\lambda+1}$, we have $\alpha=\overline{H'}\in\ASupp(\cA_{\lambda+1})\setminus\ASupp(\cA_{\lambda})$.
	
	Assume that $\lambda$ is a limit ordinal. Since $\cA_{\lambda}\subset\cA$ is the smallest Serre subcategory containing all $\cA_{\mu}$ with $\mu<\lambda$, \cref{09809844} implies that $\ASupp(\cA_{\lambda})\subset\ASpec\cA$ is the smallest open subset containing all $\ASupp(\cA_{\mu})$ with $\mu<\lambda$, which is
	\begin{equation*}
		\bigcup_{\mu<\lambda}\ASupp(\cA_{\mu})=\bigcup_{\mu<\lambda}(\ASpec\cA)_{\mu}=(\ASpec\cA)_{\lambda}
	\end{equation*}
	by the induction hypothesis.
	
	Let $\cG$ be a locally noetherian Grothendieck category and let $\cA:=\noeth\cG$. Then $\ASupp(\cG_{\lambda})$ is homeomorphic to the subsets
	\begin{equation*}
		\ASupp(\cG_{\lambda}\cap\cA)=\ASupp(\cA_{\lambda})=(\ASpec\cA)_{\lambda}
	\end{equation*}
	of $\ASpec\cA$, and the last one is homeomorphic to $(\ASpec\cG)_{\lambda}$. Since the two homeomorphisms are both induced from the homeomorphism in \cref{89431240}, we obtain the equality $\ASupp(\cG_{\lambda})=(\ASpec\cG)_{\lambda}$.
\end{proof}

\begin{Proposition}\label{63990180}
	Let $\cA$ be a locally noetherian Grothendieck category or a noetherian abelian category.
	\begin{enumerate}
		\item\label{23400278} For every object $M\in\cA$,
		\begin{equation*}
			\KGdim M=\sup\setwithcondition{\KGdim\alpha}{\alpha\in\ASupp M}.
		\end{equation*}
		\item\label{80976446} We have
		\begin{equation*}
			\KGdim\cA=\sup\setwithcondition{\KGdim\alpha}{\alpha\in\ASpec\cA}.
		\end{equation*}
	\end{enumerate}
\end{Proposition}

\begin{proof}
	This follows from \cref{09809844,60891414}.
\end{proof}

\begin{Remark}\label{95397124}
	For a locally noetherian Grothendieck category $\cG$, the Krull-Gabriel dimension of $\cG$ is not necessarily equal to the \emph{dimension} of the topological space $\ASpec\cG$ found in \cite[p.~5]{MR0463157}, which is defined to be supremum of the integers $d$ such that there exists a chain $Z_{0}\subsetneq\cdots\subsetneq Z_{d}$ of irreducible closed subsets of the topological space.
	
	We consider $\Mod^{\bbZ}k[x]$ in \cref{26229957}. It is essentially shown in \cite[Example~4.7]{MR1899866} (see also \cite[Example~3.4]{MR3351569}) that
	\begin{equation*}
		\ASpec(\Mod^{\bbZ}k[x])=\set{\overline{k[x]}}\cup\setwithcondition{\overline{S(i)}}{i\in\bbZ},
	\end{equation*}
	where $S:=k[x]/(x)$ and the atoms appearing in the right-hand side are pairwise distinct. A subset $\Phi\subset\ASpec(\Mod^{\bbZ}k[x])$ (with respect to the localizing topology) is open if and only if
	\begin{itemize}
		\item $\overline{k[x]}\notin\Phi$, or
		\item $\overline{k[x]}\in\Phi$ and there exists $i_{0}\in\bbZ$ such that $\overline{S(i)}\in\Phi$ for all $i\leq i_{0}$.
	\end{itemize}
	Indeed, if $\overline{k[x]}\notin\Phi$, then $\Phi$ consists of atoms represented by simple objects, and it hence follows from the definition of the localizing topology that $\Phi$ is open. If $\overline{k[x]}\in\Phi$ and $\Phi$ is open, then there exists a nonzero subobject $L\subset k[x]$ such that $\ASupp L\subset\Phi$. Since $L=x^{n}k[x]$ for some $n\in\bbZ_{\geq 0}$ and $\ASupp x^{n}k[x]=\setwithcondition{\overline{S(i)}}{i\leq -n}$, the subset $\Phi$ should contain all $\overline{S(i)}$ with $i\leq -n$. One can conclude that the subsets $\Phi$ with $\overline{k[x]}\in\Phi$ listed above are open by applying \cref{43980984} \cref{60221442}, as there exists $n\in\bbZ_{\geq 0}$ such that $\ASupp x^{n}k[x]\subset\Phi$ and $\Phi\setminus\set{\overline{k[x]}}$ consists of atoms represented by simple objects.
	
	Hence $\KGdim\overline{S(i)}=0$ for all $i\in\bbZ$ and $\KGdim\overline{k[x]}=1$. Consequently $\KGdim(\Mod^{\bbZ}k[x])=1$.
	
	On the other hand, every irreducible closed subsets of $\ASpec(\Mod^{\bbZ}k[x])$ (with respect to the localizing topology) consists of a single point. Thus the dimension of the topological space $\ASpec(\Mod^{\bbZ}k[x])$ in the sense of \cite[p.~5]{MR0463157} is zero.
\end{Remark}

\begin{Proposition}\label{62636912}
	Let $R$ be a commutative noetherian ring.
	\begin{enumerate}
		\item\label{62484066} $\KGdim(\Mod R)$ is equal to the Krull dimension of $R$.
		\item\label{42611408} For every $R$-module $M$, $\KGdim M$ is equal to the supremum of the lengths of chains in $\Supp M$. It is equal to the Krull dimension of $M$ if $M$ is finitely generated.
		\item\label{23133825} For every $\kp\in\Spec R$, $\KGdim\overline{R/\kp}$ is equal to the Krull dimension of $R/\kp$.
	\end{enumerate}
	In these statements, Krull-Gabriel dimensions are regarded as elements of $\bbZ_{\geq -1}\cup\set{\infty}$.
\end{Proposition}

\begin{proof}
	By \cref{30989084} \cref{30980948}, a subset of $X:=\ASpec(\Mod R)$ is open if and only if the corresponding subset of $\Spec R$ is specialization-closed via the bijection in \cref{30989084} \cref{20980444}. Thus $X_{0}$ corresponds to the set of all maximal ideals, $X_{1}$ corresponds to the set consisting of all maximal ideals and all prime ideals that are maximal among all non-maximal prime ideals, and in general, $X_{\lambda}$ corresponds to
	\begin{equation}
		\setwithcondition{\kp\in\Spec R}{\text{the Krull dimension of $R/\kp$ is at most $\lambda$}}.
	\end{equation}
	Thus \cref{23133825} in the statement holds and \cref{62484066} and \cref{42611408} follow from \cref{63990180}.
\end{proof}

\subsection{Pro-category}
\label{23989084}

When we consider the extension groups between atoms, it is useful to regard each atom as a pro-object, which is an object of the pro-category of the given abelian category. Here we recall the definition of the pro-category and its basic properties, and show that Yoneda products of extensions can be extended in terms of pro-objects (\cref{09480444}).

\begin{Definition}\label{23787594}
	Let $\cC$ be a category. Define the category $\Pro\cC$, which is called the \emph{pro-category} of $\cC$, as follows:
	\begin{enumerate}
		\item\label{41243215} Objects of $\Pro\cC$ are inverse systems (whose index sets are in bijection with small sets) in $\cC$.
		\item\label{91658269} For inverse systems $\cM=\set{M_{i}}_{i\in I}$ and $\cN=\set{N_{j}}_{j\in J}$ in $\cC$, define
		\begin{equation*}
			\Hom_{\Pro\cC}(\cM,\cN):=\varprojlim_{j\in J}\varinjlim_{i\in I}\Hom_{\cC}(M_{i},N_{j}).
		\end{equation*}
		\item\label{20980989} The composition of morphisms in $\Pro\cC$ is induced from that in $\cC$.
	\end{enumerate}
	Objects in $\Pro\cC$ are called \emph{pro-objects} in $\cC$.
\end{Definition}

\begin{Remark}\label{79591412}\leavevmode
	\begin{enumerate}
		\item\label{49804413} The pro-category $\Pro\cC$ of a category $\cC$ is the dual notion of the \emph{ind-category} $\Ind\cC$. Indeed, we have
		\begin{equation*}
			\Pro\cC=(\Ind\cC^{\op})^{\op}
		\end{equation*}
		(see the paragraph before \cite[Example~6.1.3]{MR2182076}).
		\item\label{20598092} In \cite[Definition~6.1.1]{MR2182076}, a pro-object in $\cC$ is defined to be a functor $\cC^{\op}\to\Set^{\op}$ that is isomorphic to a filtered limit of representable functors, where $\Set$ is the category of small sets. It is shown after \cite[Proposition~6.1.9]{MR2182076} that Hom-sets there can be written as in \cref{23787594} \cref{91658269}. Since every small filtered category admits a cofinal functor from a small directed set (see \cite[Theorem~1.5]{MR1294136}), there is a canonical equivalence from the pro-category defined in \cref{23787594} to the one defined in \cite[Definition~6.1.1]{MR2182076}.
	\end{enumerate}
\end{Remark}

\begin{Theorem}\label{20908444}
	Let $\cA$ be an abelian category.
	\begin{enumerate}
		\item\label{90498324} $\Pro\cA$ is an abelian category and it admits limits. Inverse limits in $\Pro\cA$ are exact.
		\item\label{09809890} The canonical functor $\cA\to\Pro\cA$, which sends each object $M\in\cA$ to the inverse system consisting of only $M$ and each morphism to the induced one, is fully faithful and exact. The essential image of the functor is closed under kernels, cokernels, and extensions.
		\item\label{04980989} If $\cA$ admits colimits, then $\Pro\cA$ admits colimits and the canonical functor $\cA\to\Pro\cA$ preserves colimits.
		\item\label{64841776} If $\cA$ is skeletally small, then $\Pro\cA$ admits a cogenerating set, and hence it is a coGrothendieck category, that is, $(\Pro\cA)^{\op}$ is a Grothendieck category.
	\end{enumerate}
\end{Theorem}

\begin{proof}
	\cite[Lemma~6.1.2, Corollary~6.1.17, Theorem~8.6.5, and Proposition~8.6.11]{MR2182076}.
\end{proof}

\begin{Remark}\label{23498098}
	We regard an abelian category $\cA$ as a full subcategory of $\Pro\cA$ via the fully faithful functor in \cref{20908444} \cref{09809890}. In \cref{15182006}, we also regard atoms in $\cA$ as objects in $\Pro\cA$.
\end{Remark}

\begin{Proposition}\label{09480444}
	Let $\cA$ be an abelian category. Let $d_{1},d_{2}\geq 0$ be integers. Then the Yoneda products
	\begin{equation*}
		\Ext_{\cA}^{d_{2}}(M,N)\times\Ext_{\cA}^{d_{1}}(L,M)\to\Ext_{\cA}^{d_{1}+d_{2}}(L,N)
	\end{equation*}
	for $L,M,N\in\cA$ induce a $\bbZ$-bilinear map
	\begin{equation*}
		\varprojlim_{k\in K}\varinjlim_{j\in J}\Ext_{\cA}^{d_{2}}(M_{j},N_{k})\times\varprojlim_{j\in J}\varinjlim_{i\in I}\Ext_{\cA}^{d_{1}}(L_{i},M_{j})\to\varprojlim_{k\in K}\varinjlim_{i\in I}\Ext_{\cA}^{d_{1}+d_{2}}(L_{i},N_{k})
	\end{equation*}
	for $\set{L_{i}}_{i\in I},\set{M_{j}}_{j\in J},\set{N_{k}}_{k\in K}\in\Pro\cA$.
\end{Proposition}

\begin{proof}
	The pro-category of the bounded derived category $\Db(\cA)$ is a $\bbZ$-linear category (see the paragraph before \cite[Proposition~8.6.2]{MR2182076}). The desired map is the composition of morphisms in $\Pro\Db(\cA)$.
\end{proof}

\begin{Remark}\label{40809844}
	Let $\cA$ be an abelian category. Then there is a canonical functor
	\begin{equation*}
		J\colon\Db(\Pro\cA)\to\Pro\Db(\cA)
	\end{equation*}
	(\cite[Theorem~15.4.3]{MR2182076}). However, this is not necessarily faithful. A counterexample is given in \cite[Exercise~15.2]{MR2182076} in terms of ind-categories: There is an abelian category $\cA$ that admits $\cM=\set{M_{i}}_{i\in I},\cN=\set{N_{j}}_{j\in J}\in\Ind\cA$ such that the map
	\begin{equation*}
		\Hom_{\Db(\Ind\cA)}(\cM,\cN[1])\to\Hom_{\Ind\Db(\cA)}(J(\cM),J(\cN[1]))
	\end{equation*}
	given by the functor $J$ is not injective. In fact, in the example given in the reference, the left-hand side is nonzero while the right-hand side is zero. Note that the map can be written as
	\begin{equation*}
		0\neq\Ext_{\Db(\Ind\cA)}^{1}(\cM,\cN)\to\varprojlim_{j\in J}\varinjlim_{i\in I}\Ext_{\cA}^{1}(M_{i},N_{j})=0.
	\end{equation*}
\end{Remark}

\section{Chasing extensions}
\label{49230598}

In this section, we prove \cref{75294329}, which is the first half part of the proof of \cref{30418873} in the introduction. For a given extension $\xi\in\Ext_{\cA}^{i}(M,N)$, where $i\geq 0$ is an integer and $M,N$ are objects in an abelian category $\cA$, we replace $M$ by a smaller subquotient $M'$ and find an extension $\xi'\in\Ext_{\cA}^{i}(M',N)$ that is related to $\xi$ via canonical maps. Repeating this process, we finally obtain an element of the naturally defined extension group $\Ext_{\cA}^{i}(\alpha,N)$ for some $\alpha\in\ASupp M$, under some noetherian assumption.

First we recall the definition of $\Ext_{\cA}^{i}(\alpha,N)$ introduced in \cite{MR3272068}.

\begin{Remark}\label{82049837}
	In order to define the extension groups between an atom and an object, we need to fix a monoform object representing the given atom. In \cite{MR3272068}, we only worked on a Grothendieck category and took $E(\alpha)$ as the representative of an atom $\alpha$. Although $E(\alpha)$ is not monoform in general, the uniformity is enough to define the extension groups.
	
	Since we will extend the definition to an arbitrary abelian category, we use the following convention.
\end{Remark}

\begin{Convention}\label{29890804}
	Let $\cA$ be an abelian category. For each $\alpha\in\ASpec\cA$, we fix a monoform object $H\in\cA$ such that $\overline{H}=\alpha$, which is referred to as the \emph{fixed representative} of $\alpha$.
\end{Convention}

\begin{Definition}\label{24980234}
	Let $\cA$ be an abelian category and $\alpha\in\ASpec\cA$. Let $H$ be the fixed representative of $\alpha$.
	\begin{enumerate}
		\item\label{24958204} For an integer $i\geq 0$, we define the functor $\Ext_{\cA}^{i}(\alpha,-)\colon\cA\to\Mod\bbZ$ that sends each object $N\in\cA$ to
		\begin{equation*}
			\Ext_{\cA}^{i}(\alpha,N):=\varinjlim_{0\neq H'\subset H}\Ext_{\cA}^{i}(H',N)
		\end{equation*}
		and each morphism in $\cA$ to the induced one. The direct limit is taken over the direct system consisting of all nonzero subobjects $H'\subset H$, together with the opposite relation of inclusion of subobjects. $\Ext^{0}_{\cA}(\alpha,-)$ is denoted by $\Hom_{\cA}(\alpha,-)$.
		\item\label{80483254} The \emph{residue field} of $\alpha$ is defined to be
		\begin{equation*}
			k(\alpha):=\Hom_{\cA}(\alpha,H).
		\end{equation*}
	\end{enumerate}
\end{Definition}

\begin{Remark}\label{20945844}
	Let $\cG$ be a Grothendieck category. The functor $\Hom_{\cG}(\alpha,-)$ and the residue field $k(\alpha)$ were defined in terms of the spectral category of $\cG$ in \cite[Definition~3.5]{MR3272068} and the functor $\Ext^{i}_{\cG}(\alpha,-)$ was introduced as the $i$-th right derived functor of $\Hom_{\cG}(\alpha,-)$, viewed as a functor $\cG\to\Mod k(\alpha)$, in \cite[Definition~4.1]{MR3272068}. It is shown in \cite[Remarks 3.6 and 4.8]{MR3272068} that those definitions are equivalent to \cref{24980234}. In particular, we have the following:
	\begin{enumerate}
		\item\label{34247366} The isomorphism class of the functor $\Ext_{\cG}^{i}(\alpha,-)\colon\cG\to\Mod\bbZ$ does not depend on the choice of the fixed representative $H$ of $\alpha$.
		\item\label{87523454} $k(\alpha)$ has a structure of a skew field, whose multiplication is induced from the composition of morphisms.
		\item\label{03598094} $\Ext_{\cG}^{i}(\alpha,N)$ has a structure of right $k(\alpha)$-module for each object $N\in\cG$, and $\Ext_{\cG}^{i}(\alpha,-)$ becomes a functor $\cG\to\Mod k(\alpha)$.
		\item\label{82893452} The functor $\Hom_{\cG}(\alpha,-)\colon\cG\to\Mod k(\alpha)$ is left exact, and $\Ext^{i}_{\cG}(\alpha,-)\colon\cG\to\Mod k(\alpha)$ is the $i$-th right derived functor of $\Hom_{\cG}(\alpha,-)$.
	\end{enumerate}
	Although $\cG$ was assumed to be locally noetherian in \cite{MR3272068}, the assumption is not necessary for any of these arguments.
\end{Remark}

The next two results show that atoms behave like noetherian objects in a locally noetherian Grothendieck category. These will be used later.

\begin{Proposition}\label{23750238}
	Let $\cG$ be a locally noetherian Grothendieck category and $\alpha\in\ASpec\cG$. Then the functor $\Hom_{\cG}(\alpha,-)\colon\cG\to\Mod k(\alpha)$ preserves direct limits.
\end{Proposition}

\begin{proof}
	Since $\cG$ be locally noetherian, we can take a nonzero noetherian subobject $H$ of the fixed representative of $\alpha$. For every direct system $\set{M_{i}}_{i\in I}$ in $\cG$, we have canonical isomorphisms
	\begin{align*}
		\Hom_{\cG}(\alpha,\varinjlim_{i\in I}M_{i})
		&\cong\varinjlim_{0\neq H'\subset H}\Hom_{\cG}(H',\varinjlim_{i\in I}M_{i})
		\cong\varinjlim_{0\neq H'\subset H}\Bigg(\varinjlim_{i\in I}\Hom_{\cG}(H',M_{i})\Bigg)\\
		&\cong\varinjlim_{i\in I}\Bigg(\varinjlim_{0\neq H'\subset H}\Hom_{\cG}(H',M_{i})\Bigg)
		\cong\varinjlim_{i\in I}\Hom_{\cG}(\alpha,M_{i}),
	\end{align*}
	where we have the second isomorphism since $H'$ is noetherian.
\end{proof}

\begin{Proposition}\label{09320089}
	Let $\cG$ be a locally noetherian Grothendieck category. Let $i\geq 0$ be an integer and $\alpha\in\ASpec\cG$. Then the functor
	\begin{equation*}
		\Ext_{\cG}^{i}(\alpha,-)\colon\cG\to\Mod k(\alpha)
	\end{equation*}
	preserves direct limits.
\end{Proposition}

\begin{proof}
	Since we have \cref{23750238}, this can be shown in a similar way to \cref{93588947}.
\end{proof}

The following is the main result in this section:

\begin{Proposition}\label{23454349}
	Let $\cA$ be an abelian category and let $0\neq\xi\in\Ext^{i}_{\cA}(M,N)$, where $i\geq 0$ is an integer, $M\in\cA$ is a noetherian object, and $N\in\cA$ is an object. Then there exist subobjects $L'\subsetneq L_{0}\subset M$ and $\eta\in\Ext_{\cA}^{i}(M/L',N)$ satisfying the following conditions:
	\begin{enumerate}
		\item\label{23985024} $L_{0}/L'$ is a monoform object. Let $\alpha\in\ASpec\cA$ be the atom represented by it.
		\item\label{34985023} For every nonzero subobject $L/L'$ of $L_{0}/L'$, the element $\eta$ is sent to
		\begin{itemize}
			\item $\xi$ in $\Ext^{i}_{\cA}(M,N)$,
			\item a nonzero element in $\Ext_{\cA}^{i}(L,N)$, and
			\item a nonzero element in $\Ext_{\cA}^{i}(\alpha,N)$
		\end{itemize}
		along the commutative diagram
		\begin{equation*}
			\begin{tikzcd}[row sep=small]
				& \Ext_{\cA}^{i}(M/L',N)\ar[dl]\ar[dr] & & \\
				\Ext_{\cA}^{i}(M,N)\ar[dr] & & \Ext_{\cA}^{i}(L/L',N)\ar[dl]\ar[dr] & \\
				& \Ext_{\cA}^{i}(L,N) & & \Ext_{\cA}^{i}(\alpha,N)\rlap{,}
			\end{tikzcd}
		\end{equation*}
		where the bottom-right map is induced from an arbitrarily fixed nonzero element of $\Hom_{\cA}(\alpha,L/L')$ and the other maps are induced from inclusions and projections.
	\end{enumerate}
\end{Proposition}

\begin{proof}
	We write
	\begin{equation*}
		F:=\Ext_{\cA}^{i}(-,N)\colon\cA^{\op}\to\Mod\bbZ\quad\text{and}\quad F(\alpha):=\Ext_{\cA}^{i}(\alpha,N).
	\end{equation*}
	The functor $F$ is half exact, that is, for every short exact sequence
	\begin{equation*}
		0\to E'\to E\to E''\to 0
	\end{equation*}
	in $\cA$, the induced sequence
	\begin{equation*}
		F(E'')\to F(E)\to F(E')
	\end{equation*}
	is exact. Moreover, $F(\alpha)$ is by definition a direct limit of $F(H)$ for various $H$. In this proof, we do not use any other property of the functor $F$.
	
	Let $L'\subset M$ be a subobject that is maximal among those satisfying the following property: There exist a subobject $L_{0}\subset M$ with $L'\subsetneq L_{0}$ and $\eta\in F(M/L')$ such that $\eta$ is sent to $\xi\in F(M)$ and a nonzero element $\zeta_{0}\in F(L_{0})$ along the above commutative diagram with $L$ replaced by $L_{0}$. Such $L'$ indeed exists because $M$ is noetherian and the zero subobject of $M$ satisfies the given property by taking $L_{0}:=M$ and $\eta:=\xi$.
	
	Let $L/L'\subset L_{0}/L'$ be a nonzero subobject. Assume that the canonical map $F(L_{0})\to F(L)$ sends $\zeta_{0}$ to zero. Then $\xi$ is sent to zero by the second map of the exact sequence
	\begin{equation*}
		F(M/L)\to F(M)\to F(L),
	\end{equation*}
	so $\xi$ is an image of some nonzero element of $F(M/L)$. This contradicts the maximality of $L'$. Hence $\zeta_{0}$ is sent to a nonzero element in $F(L)$. This means that $L$ also satisfies the requirement for $L_{0}$. Since $L_{0}/L'$ has a monoform subobject (see \cref{82094844} \cref{17457124}), we can assume that $L_{0}/L'$ itself is a monoform object by replacing $L_{0}$.
	
	By the commutativity of the diagram in the proposition, $\eta$ is also sent to a nonzero element $\zeta\in F(L/L')$. Let $\alpha$ be the atom represented by $L_{0}/L'$ and fix a nonzero element $[f]\in\Hom_{\cA}(\alpha,L/L')$ represented by $f\colon H\to L/L'$, where $H$ is a nonzero subobject of the fixed representative of $\alpha$. Assume that $\zeta$ is sent to the zero element of $F(\alpha)$ by the induced map. By the definition of $F(\alpha)$, there exists a nonzero subobject $H'\subset H$ such that $\zeta$ is sent to zero by the composition
	\begin{equation*}
		F(L/L')\to F(H)\to F(H').
	\end{equation*}
	Thus $\eta$ is also sent to zero in $F(H')$. This contradicts what we showed above for an arbitrary nonzero subobject of $L_{0}/L'$. Therefore $\zeta$ is sent to a nonzero element of $F(\alpha)$. This completes the proof.
\end{proof}

\begin{Remark}\label{16116889}
	As mentioned in the above proof, the statement of \cref{23454349} still holds after replacing all $\Ext_{\cA}^{i}(-,N)$ by an arbitrary half exact functor $F\colon\cA^{\op}\to\Mod\bbZ$ and defining $F(\alpha)$ in the same way as $\Ext_{\cA}^{i}(\alpha,N)$. This fact will be used to show \cref{49023408}.
\end{Remark}

\begin{Corollary}\label{75294329}
	Let $\cA$ be a locally noetherian Grothendieck category or a noetherian abelian category.
	\begin{enumerate}
		\item\label{43950239} Let $i\geq 0$ be an integer and $M,N\in\cA$. If $\Ext_{\cA}^{i}(M,N)\neq 0$, then there exists $\alpha\in\ASupp M$ such that $\Ext_{\cA}^{i}(\alpha,N)\neq 0$.
		\item\label{34990548} For every $N\in\cA$, we have
		\begin{equation*}
			\injdim N=\sup\setwithcondition{i\geq 0}{\textnormal{$\Ext_{\cA}^{i}(\alpha,N)\neq 0$ for some $\alpha\in\ASpec\cA$}}.
		\end{equation*}
		Therefore
		\begin{equation*}
			\gldim\cA=\sup\setwithcondition{i\geq 0}{\textnormal{$\Ext_{\cA}^{i}(\alpha,-)\neq 0$ for some $\alpha\in\ASpec\cA$}}.
		\end{equation*}
	\end{enumerate}
\end{Corollary}

\begin{proof}
	These are consequences of \cref{06236705,23454349}.
\end{proof}

\section{Virtual duals of extension groups between atoms}
\label{34982309}

In the previous section, we chased a nonzero extension in $\Ext_{\cA}^{i}(M,N)$ and found a nonzero element in $\Ext_{\cA}^{i}(\alpha,N)$ for some $\alpha\in\ASupp M$. In this section, we take the same approach for the object $N$ in the second argument. However, instead of the naturally defined extension group $\Ext_{\cA}^{i}(\alpha,\beta)$ for atoms $\alpha$ and $\beta$, we will use its \emph{virtual dual} $\varDaExt_{\cA}^{i}(\alpha,\beta)$ (see \cref{23809324}). This is because the definition of $\Ext_{\cA}^{i}(\alpha,\beta)$ involves an inverse limit, which is difficult to control, while that of the virtual dual does not.

For an abelian category $\cA$ and $\alpha\in\ASpec\cA$, we define the functor
\begin{equation*}
	D_{\alpha}:=\Hom_{k(\alpha)}(-,k(\alpha))\colon(\Mod k(\alpha))^{\op}\to\Mod (k(\alpha)^{\op}),
\end{equation*}
where $k(\alpha)^{\op}$ is the opposite skew field of $k(\alpha)$.

\begin{Definition}\label{23809324}
	Let $\cA$ be a abelian category. Let $i\geq 0$ be an integer and $\alpha,\beta\in\ASpec\cA$. Let $H$ be the fixed representative of $\beta$.
	\begin{enumerate}
		\item\label{58168196} Define
		\begin{equation*}
			\Ext_{\cA}^{i}(\alpha,\beta):=\varprojlim_{0\neq H'\subset H}\Ext_{\cA}^{i}(\alpha,H'),
		\end{equation*}
		where the direct limit is taken over all nonzero subobjects $H'\subset H$. $\Ext^{0}_{\cA}(\alpha,\beta)$ is denoted by $\Hom_{\cA}(\alpha,\beta)$.
		\item\label{23985244} Define
		\begin{equation*}
			\varDaExt_{\cA}^{i}(\alpha,\beta):=\varinjlim_{0\neq H'\subset H}\DaExt_{\cA}^{i}(\alpha,H').
		\end{equation*}
		We call it the \emph{virtual dual} of $\Ext_{\cA}^{i}(\alpha,\beta)$.
		\item\label{30598290} We say that \emph{$\varExt_{\cA}^{i}(\alpha,\beta)$ is eventually constant} (\resp \emph{eventually epic}) if there exists a nonzero subobject $H'\subset H$ such that for any nonzero subobjects $H'_{2}\subset H'_{1}\subset H'$, the canonical map
		\begin{equation*}
			\Ext_{\cA}^{i}(\alpha,H'_{2})\to\Ext_{\cA}^{i}(\alpha,H'_{1})
		\end{equation*}
		is bijective (\resp surjective).
	\end{enumerate}
\end{Definition}

\begin{Remark}\label{80372964}
	In \cref{23809324}, $\varExt_{\cA}^{i}(\alpha,\beta)$ itself is never defined. It only appears as $\varDaExt_{\cA}^{i}(\alpha,\beta)$ as in \cref{23809324} \cref{23985244} or as part of the terminologies defined in \cref{23809324} \cref{30598290}.
	
	If $\varExt_{\cA}^{i}(\alpha,\beta)$ is eventually constant (in the sense of \cref{23809324} \cref{30598290}), then $\varDaExt_{\cA}^{i}(\alpha,\beta)=\DaExt_{\cA}^{i}(\alpha,\beta)$. However, there is no guarantee that this equality holds in general.
\end{Remark}

The new notions in \cref{23809324} can easily be computed when $i=0$:

\begin{Proposition}\label{40545360}
	Let $\cA$ be an abelian category and let $\alpha,\beta\in\ASpec\cA$. Let $H$ be the fixed representative of $\beta$.
	\begin{enumerate}
		\item\label{19711748} For every nonzero subobject $H'\subset H$, the canonical map
		\begin{equation*}
			\Hom_{\cA}(\alpha,H')\to\Hom_{\cA}(\alpha,H)
		\end{equation*}
		is an isomorphism of right $k(\alpha)$-modules. 
		\item\label{35170236} The canonical map
		\begin{equation*}
			\Hom_{\cA}(\alpha,\beta)\to\Hom_{\cA}(\alpha,H)
		\end{equation*}
		is an isomorphism of right $k(\alpha)$-modules, $\varExt_{\cA}^{0}(\alpha,\beta)$ is eventually constant, and
		\begin{equation*}
			\varDaExt_{\cA}^{0}(\alpha,\beta)=\DaHom_{\cA}(\alpha,\beta).
		\end{equation*}
		\item\label{32027008} If $\alpha\neq\beta$, then $\Hom_{\cA}(\alpha,\beta)=0$. If $\alpha=\beta$, then we have a canonical isomorphism
		\begin{equation*}
			\Hom_{\cA}(\alpha,\alpha)\isoto k(\alpha)
		\end{equation*}
		of right $k(\alpha)$-modules. We identify $\Hom_{\cA}(\alpha,\alpha)$ with $k(\alpha)$ in this way.
	\end{enumerate}
\end{Proposition}

\begin{proof}
	\cref{19711748} follows from \cite[Definition~3.5]{MR3272068} (see \cref{20945844}) and \cite[Theorem~3.3 (1)]{MR3272068} since every nonzero subobject of a monoform object is an essential subobject.
	
	\cref{35170236} is an immediate consequence of \cref{19711748}.
	
	\cref{32027008} If $\alpha\neq\beta$, then $\Hom_{\cA}(\alpha,H')=0$ for all nonzero subobjects $H'\subset H$ by \cite[Proposition~3.13]{MR3272068}. If $\alpha=\beta$, then the isomorphism is identical to the one in \cref{35170236}.
\end{proof}

\begin{Remark}\label{15182006}
	Every object in an abelian category $\cA$ can be regarded as a pro-object in $\cA$ by the canonical functor $\cA\to\Pro\cA$ in \cref{20908444} \cref{09809890}. Moreover, for $\alpha\in\ASpec\cA$ with fixed representative $H$, the inverse system consisting of all nonzero subobjects of $H$ can be regarded as a pro-object in $\cA$, which will be identified with the atom $\alpha$. Hence \cref{09480444} gives natural compositions of extensions for various combinations of atoms and objects listed below. Let $i,j\geq 0$ be integers, $M,N\in\cA$, and $\alpha,\beta,\gamma\in\ASpec\cA$. We use the identification in \cref{40545360} \cref{32027008}.
	\begin{enumerate}
		\item\label{34887089} $\Ext_{\cA}^{j}(M,N)\times\Ext_{\cA}^{i}(\alpha,M)\to\Ext_{\cA}^{i+j}(\alpha,N)$.
		\item\label{01722256} $\Ext_{\cA}^{j}(\beta,N)\times\Ext_{\cA}^{i}(\alpha,\beta)\to\Ext_{\cA}^{i+j}(\alpha,N)$. Its special case
		\begin{equation*}
			\Ext_{\cA}^{i}(\alpha,N)\times k(\alpha)\to\Ext_{\cA}^{i}(\alpha,N)
		\end{equation*}
		defines a right $k(\alpha)$-module structure on $\Ext_{\cA}^{i}(\alpha,N)$, which is the same as the one in \cref{20945844}.
		\item\label{44979779} $\Ext_{\cA}^{j}(\beta,\gamma)\times\Ext_{\cA}^{i}(\alpha,\beta)\to\Ext_{\cA}^{i+j}(\alpha,\gamma)$. Its special cases
		\begin{equation*}
			\Ext_{\cA}^{i}(\alpha,\beta)\times k(\alpha)\to\Ext_{\cA}^{i}(\alpha,\beta)
		\end{equation*}
		and
		\begin{equation*}
			k(\beta)\times\Ext_{\cA}^{i}(\alpha,\beta)\to\Ext_{\cA}^{i}(\alpha,\beta)
		\end{equation*}
		define a right $(k(\beta)^{\op}\otimes_{\bbZ}k(\alpha))$-module structure on $\Ext_{\cA}^{i}(\alpha,\beta)$. Moreover, the further special case
		\begin{equation*}
			k(\alpha)\times k(\alpha)\to k(\alpha)
		\end{equation*}
		defines the structure of skew field on $k(\alpha)$, which is the same as the one in \cref{20945844}.
	\end{enumerate}
\end{Remark}

\begin{Remark}\label{10271676}
	Let $\cA$ be an abelian category. Every simple object $S\in\cA$ is monoform, and $\overline{S}\in\ASpec\cA$ is canonically isomorphic to $S$ itself as a pro-object in $\cA$. The isomorphism $\overline{S}\isoto S$ induces
	\begin{equation*}
		\Ext_{\cA}^{i}(\overline{S},-)\isofrom\Ext_{\cA}^{i}(S,-)\quad\text{and}\quad\Ext_{\cA}^{i}(\alpha,\overline{S})\isoto\Ext_{\cA}^{i}(\alpha,S)
	\end{equation*}
	for all integer $i\geq 0$ and $\alpha\in\ASpec\cA$.
\end{Remark}

The next result gives the second half part of the proof of \cref{30418873}:

\begin{Proposition}\label{49023408}
	Let $\cA$ be an abelian category and let $0\neq\xi\in\DaExt^{i}_{\cA}(\alpha,N)$, where $i\geq 0$ is an integer, $\alpha\in\ASpec\cA$, and $N\in\cA$ is a noetherian object. Then there exist subobjects $L'\subsetneq L_{0}\subset N$ and $\eta\in\DaExt_{\cA}^{i}(\alpha,N/L')$ satisfying the following conditions:
	\begin{enumerate}
		\item\label{42984984} $L_{0}/L'$ is a monoform object. Let $\beta\in\ASpec\cA$ be the atom represented by it.
		\item\label{72398454} For every nonzero subobject $L/L'$ of $L_{0}/L'$, the element $\eta$ is sent to
		\begin{itemize}
			\item $\xi$ in $\DaExt^{i}_{\cA}(\alpha,N)$,
			\item a nonzero element in $\DaExt_{\cA}^{i}(\alpha,L)$, and
			\item a nonzero element in $\varDaExt_{\cA}^{i}(\alpha,\beta)$
		\end{itemize}
		along the commutative diagram
		\begin{equation*}
			\begin{tikzcd}[row sep=small]
				& \DaExt_{\cA}^{i}(\alpha,N/L')\ar[dl]\ar[dr] & & \\
				\DaExt_{\cA}^{i}(\alpha,N)\ar[dr] & & \DaExt_{\cA}^{i}(\alpha,L/L')\ar[dl]\ar[dr] & \\
				& \DaExt_{\cA}^{i}(\alpha,L) & & \varDaExt_{\cA}^{i}(\alpha,\beta)\rlap{,}
			\end{tikzcd}
		\end{equation*}
		where the bottom-right map is induced from an arbitrarily fixed nonzero element of $\Hom_{\cA}(\beta,L/L')$ and the other maps are induced from inclusions and projections.
	\end{enumerate}
\end{Proposition}

\begin{proof}
	The functor
	\begin{equation*}
		F:=\DaExt_{\cA}^{i}(\alpha,-)\colon\cA^{\op}\to\Mod k(\alpha)^{\op}
	\end{equation*}
	is half exact. Thus, as we observed in \cref{16116889}, the proof of \cref{23454349} also works for this claim.
\end{proof}

\begin{Corollary}[\cref{30418873}]\label{24385902}
	Let $\cA$ be a locally noetherian Grothendieck category or a noetherian abelian category.
	\begin{enumerate}
		\item\label{29358034} Let $i\geq 0$ be an integer, $\alpha\in\ASpec\cA$, and $N\in\cA$. If $\Ext_{\cA}^{i}(\alpha,N)\neq 0$, then there exists $\beta\in\ASupp N$ such that $\varDaExt_{\cA}^{i}(\alpha,\beta)\neq 0$.
		\item\label{43259083} Let $i\geq 0$ be an integer and $M,N\in\cA$. If $\Ext_{\cA}^{i}(M,N)\neq 0$, then there exist $\alpha\in\ASupp M$ and $\beta\in\ASupp N$ such that $\varDaExt_{\cA}^{i}(\alpha,\beta)\neq 0$.
		\item\label{29045544} We have
		\begin{equation*}
			\gldim\cA=\sup\setwithcondition{i\geq 0}{\textnormal{$\varDaExt_{\cA}^{i}(\alpha,\beta)\neq 0$ for some $\alpha,\beta\in\ASpec\cA$}}.
		\end{equation*}
	\end{enumerate}
\end{Corollary}

\begin{proof}
	These follow from \cref{09320089,75294329,49023408}.
\end{proof}

\section{Topological properties of atom spectra}
\label{40598044}

In this section, we show some topological properties of atom spectra, which will be used in \cref{43982390}. We recall the definition of limit points, and set up some necessary notations:

\begin{Definition}\label{81958743}
	Let $X$ be a topological space and let $S\subset X$ be a subset.
	\begin{enumerate}
		\item\label{34958234} A point $x\in X$ is called a \emph{limit point} of $S$ if $x$ belongs to the closure of $S\setminus\set{x}$. The set of all limit points of $S$ is denoted by $L(S)$.
		\item\label{39485209} For an integer $i\geq 0$, define the subset $L^{i}(S)\subset X$ inductively as follows: $L^{0}(S):=S$ and $L^{i}(S):=L(L^{i-1}(S))$ if $i\geq 1$.
	\end{enumerate}
\end{Definition}

We prove some elementary results for the convenience of the reader:

\begin{Lemma}\label{43950984}
	Let $X$ be a topological space.
	\begin{enumerate}
		\item\label{30909844} For every $x\in X$, we have $L(\set{x})=\overline{\set{x}}\setminus\set{x}$.
		\item\label{30980934} For subsets $S_{1},S_{2}\subset X$, we have $L(S_{1}\cup S_{2})=L(S_{1})\cup L(S_{2})$.
		\item\label{10298558} Let $U\subset X$ be an open subset and $S\subset X$ a subset. Then
		\begin{enumerate}
			\item\label{40923183} $U\cap\overline{S}\subset\overline{U\cap S}$, and
			\item\label{29879344} $U\cap L(S)\subset L(U\cap S)$.
		\end{enumerate}
		\item\label{90671803} Let $F\subset X$ be a closed subset. Then
		\begin{equation*}
			L(F)=F\setminus\setwithtext{open points of $F$},
		\end{equation*}
		where an open point of $F$ means a point $x\in F$ such that $\set{x}$ is an open subset of the topological space $F$.
	\end{enumerate}
\end{Lemma}

\begin{proof}
	\cref{30909844} If $x\neq y\in X$, then
	\begin{equation*}
		y\in L(\set{x})\iff y\in\overline{\set{x}\setminus\set{y}}\iff y\in\overline{\set{x}}.
	\end{equation*}
	By definition, $x\notin L(\set{x})$.
	
	\cref{30980934} Let $x\in L(S_{1}\cup S_{2})$. Then
	\begin{equation*}
		x\in\overline{(S_{1}\cup S_{2})\setminus\set{x}}=\overline{(S_{1}\setminus\set{x})\cup(S_{2}\setminus\set{x})}=\overline{S_{1}\setminus\set{x}}\cup\overline{S_{2}\setminus\set{x}}.
	\end{equation*}
	Thus $x\in L(S_{1})\cup L(S_{2})$. The other inclusion is obvious.
	
	\cref{10298558} \cref{40923183} Let $x\in U\cap\overline{S}$ and assume $x\notin\overline{U\cap S}$. Then there exists an open neighborhood $V$ of $x$ such that $V\cap U\cap S=\emptyset$. Since $U\cap V$ is an open neighborhood of $x$, this means that $x\notin\overline{S}$, which is a contradiction.
	
	\cref{29879344} Let $x\in U\cap L(S)$. Then $x\in U\cap\overline{S\setminus\set{x}}$. By \cref{40923183},
	\begin{equation*}
		x\in\overline{U\cap(S\setminus\set{x})}=\overline{(U\cap S)\setminus\set{x}}.
	\end{equation*}
	This means $x\in L(U\cap S)$.
	
	\cref{90671803} Let $y\in L(F)$. Then
	\begin{equation*}
		y\in\overline{F\setminus\set{y}}\subset\overline{F}=F.
	\end{equation*}
	Thus $L(F)\subset F$. It immediately follows from the definition that a point in $F$ belongs to $L(F)$ if and only if $x$ is not an open point of $F$.
\end{proof}

\begin{Proposition}\label{93485344}
	Let $X$ be a topological space such that no point in $X$ is a limit point of its closure, that is,
	\begin{equation*}
		x\notin\overline{\overline{\set{x}}\setminus\set{x}}
	\end{equation*}
	for all $x\in X$. Let $S\subset X$ be a subset.
	\begin{enumerate}
		\item\label{40580923} $L(S)$ is a closed subset of $X$.
		\item\label{39580943} $L(\overline{S})=L(S)$.
	\end{enumerate}
\end{Proposition}

\begin{proof}
	\cref{40580923} Assume that there is an element $x\in\overline{L(S)}\setminus L(S)$. Since $x\notin\overline{S\setminus\set{x}}$, there exists an open neighborhood $U$ of $x$ such that $U\cap(S\setminus\set{x})=\emptyset$, which means $U\cap S\subset\set{x}$. Using \cref{30909844} and \cref{10298558} of \cref{43950984}, we obtain
	\begin{equation*}
		x\in U\cap\overline{L(S)}\subset\overline{U\cap L(S)}\subset\overline{L(U\cap S)}\subset\overline{L(\set{x})}=\overline{\overline{\set{x}}\setminus\set{x}}.
	\end{equation*}
	This contradicts the assumption on $X$. Therefore $\overline{L(S)}=L(S)$.
	
	\cref{39580943} Assume that there is an element $x\in L(\overline{S})\setminus L(S)$. As in \cref{40580923}, $x\notin L(S)$ implies that there exists an open neighborhood $U$ of $x$ such that $U\cap S\subset\set{x}$. By \cref{43950984} \cref{10298558}, we have
	\begin{equation*}
		x\in U\cap L(\overline{S})\subset L(U\cap\overline{S})\subset L(\overline{U\cap S})\subset L(\overline{\set{x}}).
	\end{equation*}
	This contradicts the assumption on $X$. Therefore $L(\overline{S})\subset L(S)$. The other inclusion is obvious.
\end{proof}

\begin{Remark}\label{93699633}
	Under the assumptions of \cref{93485344}, $L^{i}(S)\subset X$ is a closed subset for all $i\geq 1$. Hence we have
	\begin{equation*}
		L(S)\supset L^{2}(S)\supset\cdots
	\end{equation*}
	by \cref{43950984} \cref{90671803}.
\end{Remark}

The following result reveals a new topological property that all atom spectra possess:

\begin{Proposition}\label{29385092}
	Let $\cA$ be an abelian category. Then $\ASpec\cA$ has no point that is a limit point of its closure.
\end{Proposition}

\begin{proof}
	Let $\alpha\in\ASpec\cA$ and assume $\alpha\in\overline{\overline{\set{\alpha}}\setminus\set{\alpha}}$. Let $H\in\cA$ be a monoform object with $\overline{H}=\alpha$. Since $\ASupp H$ is an open neighborhood of $\alpha$, the set
	\begin{equation*}
		\ASupp H\cap(\overline{\set{\alpha}}\setminus\set{\alpha})
	\end{equation*}
	contains an element $\beta$. Since $\beta\in\ASupp H$ and $\beta\neq\alpha=\overline{H}$, there exists a nonzero subobject $L\subset H$ such that $\beta\in\ASupp(H/L)$. Then $\beta\in\overline{\set{\alpha}}$ implies $\alpha\in\ASupp(H/L)$. This contradicts the monoformness of $H$ (see \cite[Proposition~2.14]{MR3351569}). Therefore $\alpha\notin\overline{\overline{\set{\alpha}}\setminus\set{\alpha}}$.
\end{proof}

\cref{29385092} has a remarkable consequence. Hochster showed in \cite[Theorem~6]{MR0251026} that a topological space is homeomorphic to $\Spec R$ for some commutative ring $R$, equipped with the Zariski topology, if and only if the topological space is a spectral space defined as follows:

\begin{Definition}\label{60330213}
	A topological space $X$ is called a \emph{spectral space} if the following conditions are satisfied:
	\begin{enumerate}
		\item\label{16488501} $X$ is a Kolmogorov space.
		\item\label{20980981} $X$ is quasi-compact.
		\item\label{20989022} Finite intersections of quasi-compact open subsets of $X$ are again quasi-compact.
		\item\label{87334533} The quasi-compact open subsets of $X$ form an open basis.
		\item\label{28347523} $X$ is sober, that is, every nonempty irreducible closed subset of $X$ has a generic point.
	\end{enumerate}
\end{Definition}

For a spectral space $X$, its \emph{Hochster dual} $X^{*}$ is defined to be the topological space characterized as follows:
\begin{itemize}
	\item The underlying set of $X^{*}$ is the same as $X$.
	\item The collection of quasi-compact open subsets of $X$ is a closed basis of $X^{*}$.
\end{itemize}
It is shown in \cite[Proposition~8]{MR0251026} that the Hochster dual of a spectral space is again a spectral space. For a commutative noetherian ring $R$, all open subsets of $\Spec R$ with the Zariski topology are quasi-compact, so the Hochster dual of $\Spec R$ is homeomorphic to $\ASpec(\Mod R)$ with the localizing topology via the bijection in \cref{30989084}. Thus $\ASpec(\Mod R)$ is also a spectral space.

The next example shows that there exists a spectral space that is not homeomorphic to the atom spectrum of any abelian category.

\begin{Example}\label{43496855}
	Define a topological space $X$ by
	\begin{itemize}
		\item $X=\setwithcondition{x_{i}}{i\in\bbZ_{\geq 0}\cup\set{\infty}}$, where $x_{i}$'s are pairwise distinct.
		\item Nonempty open subsets of $X$ are of the form
		\begin{equation*}
			U_{j}:=\setwithcondition{x_{i}}{j\leq i}
		\end{equation*}
		for $j\in\bbZ_{\geq 0}$, where $\infty$ is larger than any integer. In particular, $\set{x_{\infty}}$ is \emph{not} open.
	\end{itemize}
	This space has the following properties:
	\begin{enumerate}
		\item\label{06433717} $X$ is a noetherian topological space, that is, the open subsets of $X$ satisfy the ascending chain condition. Hence every subset of $X$ is quasi-compact.
		\item\label{20909890} $X$ is a spectral space. Indeed,
		\begin{equation*}
			X\setminus U_{j}=\overline{\set{x_{j-1}}}
		\end{equation*}
		for all $0\neq j\in\bbZ_{\geq 0}$.
		\item\label{20890994} $L(\overline{\set{x_{\infty}}})=L(X)=X$ and $L(\set{x_{\infty}})=X\setminus\set{x_{\infty}}$. In particular,
		\begin{itemize}
			\item $x_{\infty}\in X$ is a limit point of its closure,
			\item $L(\set{x_{\infty}})$ is not closed, and
			\item $L(\overline{\set{x_{\infty}}})\neq L(\set{x_{\infty}})$.
		\end{itemize}
	\end{enumerate}
	All these properties immediately follow from the definition of $X$.
	
	By the property \cref{20890994} and \cref{29385092}, we conclude that $X$ is a spectral space that is not homeomorphic to the atom spectrum of any abelian category (\cref{47194013}). The property \cref{20890994} also shows that the assumption of \cref{93485344} cannot be removed.
\end{Example}

\section{Topological observation on extension groups between atoms}
\label{43982390}

The behavior of extension groups $\Ext_{\cA}^{i}(\alpha,\beta)$ between atoms are more difficult to understand than its virtual dual since it is defined as inverse limits (see \cref{28345908}). In this section, we will find a situation where $\varExt_{\cA}^{i}(\alpha,\beta)$ is eventually constant, in which case we do not have to take the inverse limit to define $\Ext_{\cA}^{i}(\alpha,\beta)$.

\begin{Definition}\label{28450923}
	Let $\cA$ be an abelian category. Let $i\geq 0$ be an integer and $\alpha\in\ASpec\cA$. Define
	\begin{equation*}
		\Epi_{i}(\alpha):=\setwithcondition{\beta\in\ASpec\cA}{\textnormal{$\varExt_{\cA}^{i}(\alpha,\beta)$ is eventually epic and $\varDaExt_{\cA}^{i}(\alpha,\beta)\neq 0$}}
	\end{equation*}
	and
	\begin{equation*}
		\Const_{i}(\alpha):=\setwithcondition{\beta\in\ASpec\cA}{\textnormal{$\varExt_{\cA}^{i}(\alpha,\beta)$ is eventually constant and $\Ext_{\cA}^{i}(\alpha,\beta)\neq 0$}}.
	\end{equation*}
\end{Definition}

By using a similar approach to \cref{23454349}, we obtain the following result:

\begin{Proposition}\label{92347589}
	Let $\cA$ be an abelian category. Let $i\geq 0$ be an integer, $\alpha\in\ASpec\cA$, and $N\in\cA$ a noetherian object. If $\Ext_{\cA}^{i}(\alpha,N)\neq 0$, then there exists $\beta\in\ASupp N$ such that $\varExt_{\cA}^{i}(\alpha,\beta)$ is eventually epic and $\varDaExt_{\cA}^{i}(\alpha,\beta)\neq 0$, that is,
	\begin{equation*}
		\ASupp N\cap\Epi_{i}(\alpha)\neq\emptyset.
	\end{equation*}
\end{Proposition}

\begin{proof}
	Let $L'\subset N$ be a subobject that is maximal with respect to the property that there exists a subobject $L/L'$ of $N/L'$ such that $\Ext_{\cA}^{i}(\alpha,L/L')\neq 0$. Such $L'$ exists because $N$ is noetherian and the zero subobject of $N$ satisfies the given property by taking $L:=N$. Let $L_{2}/L'\subset L_{1}/L'$ be nonzero subobjects of $L/L'$. Since we have an exact sequence
	\begin{equation*}
		\Ext_{\cA}^{i}(\alpha,L_{2}/L')\to\Ext_{\cA}^{i}(\alpha,L_{1}/L')\to\Ext_{\cA}^{i}(\alpha,L_{1}/L_{2})
	\end{equation*}
	and the last term is zero by the maximality of $L'$, the first map is surjective. In particular, by setting $L_{1}=L$, we have $\Ext_{\cA}^{i}(\alpha,L_{2}/L')\neq 0$, which means that an arbitrary nonzero subobject $L_{2}/L'$ of $L/L'$ also satisfies the requirement of $L/L'$. Since $L/L'$ has a monoform subobject, we can assume that $L/L'$ itself is a monoform object. Let $\beta\in\ASpec\cA$ be the atom represented by it. Then $\beta\in\ASupp N$, and the surjectivity mentioned above shows that $\beta\in\Epi_{i}(\alpha)$.
\end{proof}

\begin{Remark}\label{28345908}
	In general, a direct system of nonzero right modules over a ring whose structure morphisms are injective has nonzero direct limit. (Indeed, it is easy to see that any nonzero element of any module in the direct system gives a nonzero element of the direct limit.) Hence, if $\varExt_{\cA}^{i}(\alpha,\beta)$ is eventually epic in the setting of \cref{28450923}, then $\varDaExt_{\cA}^{i}(\alpha,\beta)\neq 0$ unless $\varExt_{\cA}^{i}(\alpha,\beta)$ is eventually constant and the constant value is zero.
	
	However, we do not know whether $\beta\in\Epi_{i}(\alpha)$ implies that $\Ext_{\cA}^{i}(\alpha,\beta)\neq 0$. It is shown in \cite[section~3]{MR0061086} that there is an example of an inverse system of nonzero vector spaces over a field such that all structure morphisms are surjective but the inverse limit is zero.
\end{Remark}

For this reason, we focus on the case where $\varExt_{\cA}^{i}(\alpha,\beta)$ is eventually constant. The following result is the key to the subsequent observations:

\begin{Theorem}\label{84509328}
	Let $\cA$ be an abelian category that admits a generating set consisting of noetherian objects. Let $\alpha\in\ASpec\cA$.
	\begin{enumerate}
		\item\label{30980984} $\Epi_{0}(\alpha)=\Const_{0}(\alpha)=\set{\alpha}$.
		\item\label{28389880} For every integer $i\geq 1$, we have $\Epi_{i}(\alpha)\subset\Const_{i}(\alpha)\cup L(\Epi_{i-1}(\alpha))$.
		\item\label{47932894} Consequently, we have
		\begin{equation*}
			\Epi_{i}(\alpha)\subset\bigcup_{j=0}^{i}L^{j}(\Const_{i-j}(\alpha))
		\end{equation*}
		for all integers $i\geq 0$.
	\end{enumerate}
\end{Theorem}

\begin{proof}
	\cref{30980984} This follows from \cref{40545360}.
	
	\cref{28389880} Let $\beta\in\Epi_{i}(\alpha)\setminus\Const_{i}(\alpha)$ and let $\Phi\subset\ASpec\cA$ be an open neighborhood of $\beta$. Then there exists a noetherian monoform object $H\in\cA$ such that $\overline{H}=\beta$ and $\ASupp H\subset\Phi$ (\cref{43980984} \cref{60221442}). If $\varExt_{\cA}^{i}(\alpha,\beta)$ is eventually constant, then $\varDaExt_{\cA}^{i}(\alpha,\beta)=\DaExt_{\cA}^{i}(\alpha,\beta)$ as observed in \cref{80372964}, and thus the assumption $\beta\in\Epi_{i}(\alpha)$ implies that $\Ext_{\cA}^{i}(\alpha,\beta)\neq 0$. But this contradicts the assumption $\beta\notin\Const_{i}(\alpha)$. So $\varExt_{\cA}^{i}(\alpha,\beta)$ is not eventually constant. Therefore there exist nonzero subobjects $H_{2}\subset H_{1}\subset H$ such that the canonical map $\Ext_{\cA}^{i}(\alpha,H_{2})\to\Ext_{\cA}^{i}(\alpha,H_{1})$ is not bijective but surjective (because $\varExt_{\cA}^{i}(\alpha,\beta)$ is eventually epic). Since we have an exact sequence
	\begin{equation*}
		\Ext_{\cA}^{i-1}(\alpha,H_{1}/H_{2})\to\Ext_{\cA}^{i}(\alpha,H_{2})\to\Ext_{\cA}^{i}(\alpha,H_{1}),
	\end{equation*}
	it follows that $\Ext_{\cA}^{i-1}(\alpha,H_{1}/H_{2})\neq 0$. Hence \cref{92347589} implies that
	\begin{equation*}
		\emptyset\neq\ASupp(H_{1}/H_{2})\cap\Epi_{i-1}(\alpha)\subset\ASupp H\cap\Epi_{i-1}(\alpha)\subset\Phi\cap\Epi_{i-1}(\alpha).
	\end{equation*}
	Since $H$ is monoform, $\ASupp(H_{1}/H_{2})$ does not contain $\beta$ (see \cite[Proposition~2.14]{MR3351569}). Therefore $\beta\in L(\Epi_{i-1}(\alpha))$.
	
	\cref{47932894} Using \cref{28389880} and \cref{43950984} \cref{30980934} repeatedly, we have
	\begin{align*}
		\Epi_{i}(\alpha)
		&\subset\Const_{i}(\alpha)\cup L(\Epi_{i-1}(\alpha))\\
		&\subset\Const_{i}(\alpha)\cup L(\Const_{i-1}(\alpha))\cup L^{2}(\Epi_{i-2}(\alpha))\\
		&\subset\cdots\subset\Const_{i}(\alpha)\cup L(\Const_{i-1}(\alpha))\cup\cdots\cup L^{i}(\Epi_{0}(\alpha)).
	\end{align*}
	Thus the claim follows from \cref{30980984}.
\end{proof}

\begin{Corollary}\label{23465982}
	Let $\cA$ be an abelian category that admits a generating set consisting of noetherian objects. Let $i\geq 0$ be an integer, $\alpha\in\ASpec\cA$, and $N\in\cA$ a object. If $\Ext_{\cA}^{i}(\alpha,N)\neq 0$, then there exists $\beta\in\ASupp N$ such that
	\begin{equation*}
		\beta\in\bigcup_{j=0}^{i}L^{j}(\Const_{i-j}(\alpha)).
	\end{equation*}
\end{Corollary}

\begin{proof}
	By \cref{09320089}, we can assume that $N$ is a noetherian object. Thus the claim follows from \cref{92347589,84509328}.
\end{proof}

In the case where $i=1$, the conclusion of \cref{23465982} becomes significantly simpler:

\begin{Corollary}\label{20707733}
	Let $\cA$ be an abelian category that admits a generating set consisting of noetherian objects. Let $\alpha\in\ASpec\cA$ and $N\in\cA$ an object. If $\Ext_{\cA}^{1}(\alpha,N)\neq 0$, then either $\ASupp N\cap\Const_{1}(\alpha)\neq\emptyset$ or $\alpha\in\ASupp N$.
\end{Corollary}

\begin{proof}
	By \cref{23465982}, we have
	\begin{equation*}
		\ASupp N\cap\Const_{1}(\alpha)\neq\emptyset\quad\text{or}\quad\ASupp N\cap L(\Const_{0}(\alpha))\neq\emptyset.
	\end{equation*}
	If the first assertion holds, then the result follows. So assume that the latter assertion holds. By \cref{84509328} \cref{30980984} and \cref{43950984} \cref{30909844},
	\begin{equation*}
		L(\Const_{0}(\alpha))=L(\set{\alpha})=\overline{\set{\alpha}}\setminus\set{\alpha}.
	\end{equation*}
	Hence $\ASupp N\cap\overline{\set{\alpha}}\neq\emptyset$. Since $\ASupp N$ is an open subset of $\ASpec\cA$, we have $\alpha\in\ASupp N$.
\end{proof}

In the rest of this section, we will interpret \cref{23465982} using the Krull-Gabriel dimension of an abelian category, and evaluate the difference of the following two invariants defined for atoms:

\begin{Definition}\label{96983306}
	Let $\cA$ be an abelian category and let $\alpha\in\ASpec\cA$.
	\begin{enumerate}
		\item\label{09669227} Define the \emph{projective dimension} of $\alpha$ to be
		\begin{equation*}
			\projdim\alpha:=\sup\setwithcondition{i\geq 0}{\Ext_{\cA}^{i}(\alpha,-)\neq 0},
		\end{equation*}
		where $\Ext_{\cA}^{i}(\alpha,-)$ is regarded as a functor $\cA\to\Mod k(\alpha)$.
		\item\label{42093908} Define
		\begin{equation*}
			\cprojdim\alpha:=\sup\setwithcondition{i\geq 0}{\Const_{i}(\alpha)\neq\emptyset}.
		\end{equation*}
	\end{enumerate}
\end{Definition}

\begin{Lemma}\label{17460616}
	Let $\cA$ be a locally noetherian Grothendieck category or a noetherian abelian category. Let $\Psi\subset\ASpec\cA$ be a subset. If all atoms in $\Psi$ have Krull-Gabriel dimensions at least $i$, where $i\geq 0$ is an integer, then all atoms in $L(\Psi)$ have Krull-Gabriel dimensions at least $i+1$.
\end{Lemma}

\begin{proof}
	Let $X:=\ASpec\cA$ and use the notation in \cref{38948897}. Since $X_{i}$ is the set of all atoms whose Krull-Gabriel dimensions are at most $i$, it suffices to prove that $L(X\setminus X_{i-1})\subset X\setminus X_{i}$. By definition, $X_{i-1}\subset X$ is an open subset. Hence by \cref{43950984} \cref{90671803},
	\begin{align*}
		L(X\setminus X_{i-1})
		&=(X\setminus X_{i-1})\setminus\setwithtext{open points of $X\setminus X_{i-1}$}\\
		&=(X\setminus X_{i-1})\setminus(X_{i}\setminus X_{i-1})\\
		&=X\setminus X_{i}.\qedhere
	\end{align*}
\end{proof}

\begin{Theorem}[\cref{71731514}]\label{26816538}
	Let $\cA$ be a locally noetherian Grothendieck category or a noetherian abelian category. For every $\alpha\in\ASpec\cA$, we have
	\begin{equation*}
		\projdim\alpha\leq\cprojdim\alpha+\KGdim\cA.
	\end{equation*}
\end{Theorem}

\begin{proof}
	Let $i$ be an integer such that $\Ext_{\cA}^{i}(\alpha,N)\neq 0$ for some object $N\in\cA$. It suffices to show that $i$ is less than or equal to the right-hand side of the formula. By \cref{23465982}, $L^{j}(\Const_{i-j}(\alpha))\neq\emptyset$ for some integer $0\leq j\leq i$. Set $d:=\KGdim\cA$. Then by \cref{17460616}, there exists $\beta\in\Const_{i-j}(\alpha)$ with $\KGdim\beta\leq d-j$. In particular, $\cprojdim\alpha\geq i-j$. Hence
	\begin{equation*}
		\cprojdim\alpha+\KGdim\cA\geq (i-j)+d\geq i+\KGdim\beta\geq i.\qedhere
	\end{equation*}
\end{proof}

\section{Localizing subcategories closed under injective envelopes}
\label{94285024}

Although extension groups between atoms are not easy to control in general, the first extension groups have quite a nice property described in \cref{20707733}. This allows us to determine which localizing subcategories of a locally noetherian Grothendieck category are closed under injective envelopes, in terms of extension groups between atoms.

For every object $M$ in a Grothendieck category $\cG$, its injective envelope $E(M)$ is an essential extension of $M$. On the other hand, every essential extension of $M$ is isomorphic to some subobject of $E(M)$. Thus, a localizing subcategory of $\cG$ is closed under injective envelopes if and only if it is closed under essential extensions. We state the next result using essential extensions since it also makes sense for noetherian abelian categories.

\begin{Lemma}\label{30980494}
	Let $\cA$ be a locally noetherian Grothendieck category or a noetherian abelian category. Let $\Phi\subset\ASpec\cA$ be an open subset. Then the following conditions are equivalent:
	\begin{enumerate}
		\item\label{16026888} $\ASupp^{-1}\Phi$ is closed under essential extensions.
		\item\label{55714305} $\Ext_{\cA}^{1}(\alpha,N)=0$ for all $\alpha\in\ASpec\cA\setminus\Phi$ and $N\in\ASupp^{-1}\Phi$.
		\item\label{03980944} $\varDaExt_{\cA}^{1}(\alpha,\beta)=0$ for all $\alpha\in\ASpec\cA\setminus\Phi$ and $\beta\in\Phi$.
		\item\label{39809230} $\Ext_{\cA}^{1}(\alpha,\beta)=0$ for all $\alpha\in\ASpec\cA\setminus\Phi$ and $\beta\in\Phi$.
		\item\label{30989044} $\Const_{1}(\alpha)\cap\Phi=\emptyset$ for all $\alpha\in\ASpec\cA\setminus\Phi$.
	\end{enumerate}
\end{Lemma}

\begin{proof}
	\cref{16026888}$\Rightarrow$\cref{55714305}: Assume that $\Ext_{\cA}^{1}(\alpha,N)\neq 0$ for some $\alpha\in\ASpec\cA\setminus\Phi$ and $N\in\ASupp^{-1}\Phi$. Let
	\begin{equation*}
		0\to N\to E\to H\to 0
	\end{equation*}
	be the extension corresponding to an element $\xi\in\Ext_{\cA}^{1}(H,N)$ that defines a nonzero element of $\Ext_{\cA}^{1}(\alpha,N)$, where $H\in\cA$ is a nonzero subobject of the fixed representative of $\alpha$. We regard $N$ as a subobject of $E$. If $N\subset E$ is not essential, then there exists a nonzero subobject $N'\subset E$ such that $N\cap N'=0$. Let $E':=N\oplus N'\subset E$ and $H':=E'/N\cong N'$. The element $\xi\in\Ext_{\cA}^{1}(H,N)$ is sent to the element of $\Ext_{\cA}^{1}(H',N)$ corresponding to the split exact sequence
	\begin{equation*}
		0\to N\to E'\to H'\to 0.
	\end{equation*}
	This contradicts to the fact that $\xi$ defines a nonzero element of $\Ext_{\cA}^{1}(\alpha,N)$. Hence $N\subset E$ is an essential subobject. Since $\alpha\notin\Phi$, the quotient $H=E/N$ does not belong to $\ASupp^{-1}\Phi$. Thus $\ASupp^{-1}\Phi$ is not closed under essential extensions. This is a contradiction. So $\Ext_{\cA}^{1}(\alpha,N)=0$ as desired.
	
	\cref{55714305}$\Rightarrow$\cref{03980944},\cref{39809230}: Assume \cref{55714305} and let $\alpha\in\ASpec\cA\setminus\Phi$ and $\beta\in\Phi$. Since $\Phi$ is open, there exists a monoform object $H\in\cA$ such that $\overline{H}=\beta$ and $\ASupp H\subset\Phi$ (\cref{43980984} \cref{60221442}). For every nonzero subobject $H'\subset H$, we have $\Ext_{\cA}^{1}(\alpha,H')=0$ by \cref{55714305}. Hence $\varDaExt_{\cA}^{1}(\alpha,\beta)=0$ and $\Ext_{\cA}^{1}(\alpha,\beta)=0$.
	
	\cref{03980944}$\Rightarrow$\cref{30989044} and \cref{39809230}$\Rightarrow$\cref{30989044} are obvious.
	
	\cref{30989044}$\Rightarrow$\cref{16026888}: Assume that $\ASupp^{-1}\Phi$ is not closed under essential extensions. Then there exist an object $N\in\cA$ and an essential subobject $N\subset E$ such that $N\in\ASupp^{-1}\Phi$ and $E\notin\ASupp^{-1}\Phi$. By replacing $N$, we can assume that $N$ is the largest subobject of $E$ among those belonging to $\ASupp^{-1}\Phi$. Indeed, if $\cA$ is a locally noetherian Grothendieck category, then $\ASupp^{-1}\Phi$ is a localizing subcategory, and in particular, it is closed under arbitrary sums of subobjects. If $\cA$ is a noetherian abelian category, then $\ASupp^{-1}\Phi$ is closed under finite sums of subobjects and the noetherianity of $E$ ensures the existence of such $N$. By replacing $E$ by its suitable subobject containing $N$, we can also assume that $H:=E/N$ is monoform. Since $\ASupp^{-1}\Phi$ is closed under extensions, no nonzero subobject of $H$ belongs to $\ASupp^{-1}\Phi$ by the maximality of $N$. Thus $\alpha:=\overline{H}\notin\Phi$. For every nonzero subobject $H'=E'/N$ of $H=E/N$, the element of $\Ext_{\cA}^{1}(H,N)$ corresponding to the short exact sequence
	\begin{equation*}
		0\to N\into E\onto H\to 0
	\end{equation*}
	is sent to the element of $\Ext_{\cA}^{1}(H',N)$ corresponding to the nonsplit short exact sequence
	\begin{equation*}
		0\to N\into E'\onto H'\to 0
	\end{equation*}
	by the map induced from the inclusion $H'\into H$. This means that these elements define a nonzero element of $\Ext_{\cA}^{1}(\alpha,N)$. By \cref{20707733}, we have that either $\ASupp N\cap\Const_{1}(\alpha)\neq\emptyset$ or $\alpha\in\ASupp N$. Since $\alpha\notin\Phi$ and $\ASupp N\subset\Phi$, the latter one does not hold, and hence $\Const_{1}(\alpha)\cap\Phi\neq\emptyset$. This is a contradiction. Therefore $\ASupp^{-1}\Phi$ is closed under essential extensions as desired.
\end{proof}

\begin{Theorem}[\cref{92828160}]\label{07240397}\leavevmode
	\begin{enumerate}
		\item\label{28349054} Let $\cG$ be a locally noetherian Grothendieck category. Then the bijection in \cref{09809844} \cref{93450984} induces an order-preserving bijective correspondence between
		\begin{itemize}
			\item localizing subcategories of $\cG$ that are closed under injective envelopes, and
			\item open subsets $\Phi$ of $\ASpec\cG$ with $\Ext_{\cG}^{1}(\alpha,\beta)=0$ for all $\alpha\in\ASpec\cG\setminus\Phi$ and $\beta\in\Phi$.
		\end{itemize}
		\item\label{28734504} Let $\cA$ be a noetherian abelian category. Then the bijection in \cref{09809844} \cref{42375944} induces an order-preserving bijective correspondence between
		\begin{itemize}
			\item Serre subcategories of $\cA$ that are closed under essential extensions, and
			\item open subsets $\Phi$ of $\ASpec\cA$ with $\Ext_{\cA}^{1}(\alpha,\beta)=0$ for all $\alpha\in\ASpec\cA\setminus\Phi$ and $\beta\in\Phi$.
		\end{itemize}
	\end{enumerate}
\end{Theorem}

\begin{proof}
	These follow from \cref{09809844,30980494}.
\end{proof}

\begin{Example}\label{23958094}
	Consider $\cG:=\Mod^{\bbZ}k[x]$ as in \cref{95397124}. We compute $\Ext_{\cG}^{1}(\alpha,\beta)$ for all $\alpha,\beta\in\ASpec\cG$.
	
	Let $\alpha=\overline{k[x]}$. Since every nonzero subobject $H\subset k[x]$ is isomorphic to $k[x](j)$ for some integer $j\leq 0$, it is projective, and hence $\Ext_{\cG}^{1}(H,-)=0$. Thus, for all $\beta\in\ASpec\cG$, $\varExt_{\cG}^{1}(\overline{k[x]},\beta)$ is eventually constant and $\Ext_{\cG}^{1}(\overline{k[x]},\beta)=0$.
	
	Let $\alpha=\overline{S(i)}$ for some $i\in\bbZ$ and $\beta=\overline{k[x]}$. For all integers $j\leq i-2$, we have $\Ext_{\cG}^{1}(S(i),k[x](j))=0$. Indeed, if we have an extension
	\begin{equation*}
		0\to k[x](j)\to E\to S(i)\to 0,
	\end{equation*}
	then $E_{-i}=k=E_{-j}$ and $E_{-i+1}=0$. Since $k[x]$ is generated in degree $1$ as a $k$-algebra, the extension splits. Therefore $\varExt_{\cG}^{1}(\overline{S(i)},\overline{k[x]})$ is eventually constant and $\Ext_{\cG}^{1}(\overline{S(i)},\overline{k[x]})=0$.
	
	Let $\alpha=\overline{S(i)}$ and $\beta=\overline{S(j)}$ for some $i,j\in\bbZ$. Since these atoms are represented by simple objects, $\varExt_{\cG}^{1}(\overline{S(i)},\overline{S(j)})$ is eventually constant, and
	\begin{equation*}
		\Ext_{\cG}^{1}(\overline{S(i)},\overline{S(j)})\isoto\Ext_{\cG}^{1}(S(i),S(j))\cong
		\begin{cases}
			k & \text{if $j=i-1$,} \\
			0 & \text{otherwise.}
		\end{cases}
	\end{equation*}
	
	To summarize, $\varExt_{\cG}^{1}(\alpha,\beta)$ are eventually constant for all $\alpha,\beta\in\ASpec\cG$, and
	\begin{equation*}
		\Ext_{\cG}^{1}(\alpha,\beta)\cong
		\begin{cases}
			k & \text{if $\alpha=\overline{S(i)}$ and $\beta=\overline{S(i-1)}$ for some $i\in\bbZ$,} \\
			0 & \text{otherwise.}
		\end{cases}
	\end{equation*}
	We described all open subsets of $\ASpec\cG$ in \cref{95397124}. Among them, those satisfying the conditions in \cref{30980494} are
	\begin{itemize}
		\item the empty set,
		\item $\setwithcondition{\overline{S(i)}}{i_{0}\leq i}$, where $i_{0}\in\bbZ$,
		\item $\setwithcondition{\overline{S(i)}}{i\in\bbZ}$, and
		\item $\ASpec\cG$ itself.
	\end{itemize}
	Therefore, by \cref{07240397}, all localizing subcategories of $\Mod^{\bbZ}k[x]$ closed under injective envelopes are
	\begin{itemize}
		\item the zero subcategory,
		\item $\setwithcondition{M\in\Mod^{\bbZ}k[x]}{\textnormal{$M_{i}=0$ for all $i>i_{0}$}}$, where $i_{0}\in\bbZ$,
		\item the full subcategory of $\Mod^{\bbZ}k[x]$ consisting of all torsion modules, and
		\item $\Mod^{\bbZ}k[x]$ itself,
	\end{itemize}
	and all Serre subcategories of $\mod^{\bbZ}k[x]$ closed under essential extensions are
	\begin{itemize}
		\item the zero subcategory,
		\item $\setwithcondition{M\in\mod^{\bbZ}k[x]}{\textnormal{$M_{i}=0$ for all $i>i_{0}$}}$, where $i_{0}\in\bbZ$,
		\item the full subcategory of $\mod^{\bbZ}k[x]$ consisting of all modules of finite length, and
		\item $\mod^{\bbZ}k[x]$ itself.
	\end{itemize}
\end{Example}

\begin{Remark}\label{55812411}
	Let $\Lambda$ be a finite dimensional algebra over an algebraically closed field $k$. Since $\Lambda$ is a right (and left) artinian ring, all atoms in $\Mod\Lambda$ are represented by simple objects and we have a canonical bijection between $\ASpec(\Mod\Lambda)$ and the set of isomorphism classes of simple right $\Lambda$-modules (\cite[Proposition~8.2]{MR2964615}).
	
	It is known that $\Mod\Lambda$ is equivalent to $\Mod\Lambda'$ as a $k$-linear abelian category for some \emph{basic} finite dimensional $k$-algebra $\Lambda'$ (\cite[Corollary~I.6.10]{MR2197389}). Assume that $\Lambda$ itself is a basic algebra. Then we have an isomorphism $\Lambda\cong kQ/I$ of $k$-algebras, where $Q$ is the Gabriel quiver of $\Lambda$ and $I$ is an admissible ideal of $kQ$ (\cite[Theorem~II.3.7]{MR2197389}; the non-connected case can easily be reduced to the connected case). There is a canonical bijection between the set of isomorphism classes of simple right $\Lambda$-modules and the set of vertices of $Q$ (\cite[Lemma~III.2.1 (b)]{MR2197389}). Let $S(i)$ and $S(j)$ be the simple right $\Lambda$-modules corresponding to vertices $i$ and $j$ in $Q$, respectively. By \cref{10271676}, the first extension group is
	\begin{equation*}
		\Ext_{\Lambda}^{1}(\overline{S(i)},\overline{S(j)})\isoto\Ext_{\Lambda}^{1}(S(i),S(j))
	\end{equation*}
	and its dimension over $k$ is equal to the number of arrows $i\to j$ in the Gabriel quiver $Q$ (\cite[Lemma~III.2.12 (b)]{MR2197389}).
\end{Remark}

\section{The case of noetherian algebras}
\label{48235092}

In this section, we will describe the extension groups between atoms for noetherian algebras.

Let $R$ be a commutative noetherian ring. A \emph{noetherian $R$-algebra} is a ring $\Lambda$ whose center contains $R$ as a subring such that $\Lambda$ is finitely generated as an $R$-module.

For a noetherian $R$-algebra $\Lambda$, the set of prime (two-sided) ideals of $\Lambda$ is denoted by $\Spec\Lambda$. Recall that for each $P\in\Spec\Lambda$, the intersection $P\cap R$ is a prime ideal of $R$. For every $\kp\in\Spec R$, the localization $\Lambda_{\kp}$ is a noetherian $R_{\kp}$-algebra.

For $M\in\Mod\Lambda$, we denote by $\Ass_{\Lambda}M$ the set of associated primes of $M$, that is, $P\in\Spec\Lambda$ belongs to $\Ass_{\Lambda}M$ if and only if there exists a nonzero $\Lambda$-submodule $L\subset M$ such that $\Ann_{\Lambda}L'=P$ for all nonzero $\Lambda$-submodules $L'\subset L$.

The atom spectrum of $\Mod\Lambda$ and its structure can be described in terms of the prime ideals:

\begin{Proposition}\label{63482658}
	Let $\Lambda$ be a noetherian $R$-algebra.
	\begin{enumerate}
		\item\label{13636275} There is a bijection
		\begin{equation*}
			\Spec\Lambda\isoto\ASpec(\Mod\Lambda)
		\end{equation*}
		that sends $P\in\Spec\Lambda$ to $\widetilde{P}\in\ASpec(\Mod\Lambda)$ characterized by $\AAss(\Lambda/P)=\set{\widetilde{P}}$.
		\item\label{14214126} For every $M\in\Mod\Lambda$, the bijection in \cref{13636275} induces a bijection $\Ass_{\Lambda}M\isoto\AAss M$.
		\item\label{31747619} For every $P\in\Spec\Lambda$, the bijection in \cref{13636275} induces a bijection
		\begin{equation*}
			\setwithcondition{Q\in\Spec\Lambda}{P\subset Q}\isoto\ASupp(\Lambda/P),
		\end{equation*}
		and $\ASupp(\Lambda/P)$ is the smallest open subset of $\ASpec(\Mod\Lambda)$ containing $\widetilde{P}$.
		\item\label{69982469} The bijection in \cref{13636275} is an isomorphism of partially ordered sets:
		\begin{equation*}
			(\Spec\Lambda,\subset)\isoto(\ASpec(\Mod\Lambda),\leq).
		\end{equation*}
	\end{enumerate}
\end{Proposition}

\begin{proof}
	\cref{13636275} and \cref{14214126} are shown in \cite[Theorem~7.2]{MR3272068}.
	
	\cref{31747619} Let $Q\in\Spec\Lambda$ and assume $P\subset Q$. Since we have the canonical surjection $\Lambda/P\onto\Lambda/Q$,
	\begin{equation*}
		\ASupp\frac{\Lambda}{P}\supset\ASupp\frac{\Lambda}{Q}\supset\AAss\frac{\Lambda}{Q}=\set{\widetilde{Q}}.
	\end{equation*}
	Thus $\widetilde{Q}\in\ASupp(\Lambda/P)$.
	
	Conversely, if $\widetilde{Q}\in\ASupp(\Lambda/P)$, then there exists a monoform subquotient $H$ of $\Lambda/P$ such that $\overline{H}=\widetilde{Q}$. Since $\widetilde{Q}\in\AAss H$, \cref{14214126} implies $Q\in\Ass_{\Lambda}H$. By \cite[Lemma~2.5.1]{MR1898632}, there exists a $\Lambda$-monomorphism $\Lambda/Q\to H^{\oplus n}$ for some integer $n\geq 1$. Hence $\Lambda/Q$ is isomorphic to a subquotient of $(\Lambda/P)^{\oplus n}$, and it implies that $\Lambda/Q$ is annihilated by $P$. Therefore $P\subset Q$.
	
	If $\Phi\subset\ASpec(\Mod\Lambda)$ is an open subset containing $\widetilde{P}$, then there exists a monoform right $\Lambda$-module $H$ such that $\overline{H}=\widetilde{P}$ and $\ASupp H\subset\Phi$ (\cref{43980984} \cref{60221442}). The same argument as above using \cite[Lemma~2.5.1]{MR1898632} shows that there exists a $\Lambda$-monomorphism $\Lambda/P\to H^{\oplus m}$ for some integer $m\geq 1$. Hence
	\begin{equation*}
		\ASupp\frac{\Lambda}{P}\subset\ASupp H^{\oplus m}=\ASupp H\subset\Phi.
	\end{equation*}
	This completes the proof.
	
	\cref{69982469} This follows from \cref{31747619}.
\end{proof}

We define the modules $S(P)$ and recall some basic properties, which will be used to describe the extension groups between atoms. Most of the properties are essentially observed in \cite[sections~2.4 and~2.5]{MR1898632}.

\begin{Definition}\label{14164284}
	Let $\Lambda$ be a noetherian $R$-algebra and let $P\in\Spec\Lambda$.
	\begin{enumerate}
		\item\label{37391998} Define $I(P)\in\Mod\Lambda$ to be the injective envelope $E(\widetilde{P})$ of $\widetilde{P}\in\ASpec(\Mod\Lambda)$.
		\item\label{46990266} Define a $\Lambda$-submodule $S(P)\subset I(P)$ by
		\begin{equation*}
			S(P):=\setwithcondition{x\in I(P)}{xP=0}.
		\end{equation*}
	\end{enumerate}
\end{Definition}

\begin{Proposition}\label{63542092}
	Let $\Lambda$ be a noetherian $R$-algebra. Let $P\in\Spec\Lambda$ and $\kp:=P\cap R$.
	\begin{enumerate}
		\item\label{20619540} The canonical maps $I(P)\to I(P)_{\kp}$ and $S(P)\to S(P)_{\kp}$ are isomorphisms of right $\Lambda$-modules. Thus $I(P)$ and $S(P)$ can also be regarded as right $\Lambda_{\kp}$-modules.
		\item\label{20391903} $S(P)$ is a monoform $\Lambda$-submodule of $I(P)$.
		\item\label{13184618} $S(P)$ is the unique simple $\Lambda_{\kp}$-submodule of $I(P)$.
		\item\label{75132773} For every nonzero $\Lambda$-submodule $H\subset S(P)$, we have $H_{\kp}=S(P)$.
	\end{enumerate}
\end{Proposition}

\begin{proof}
	\cite[Lemma~7.9]{MR3272068} shows that $S(P)$ defined in \cref{14164284} is isomorphic to $S(P)$ defined in the paragraph before \cite[Theorem~7.6]{MR3272068}, which is a simple right $\Lambda_{\kp}$-module. \cref{20619540} and \cref{13184618} are shown in the proof of \cite[Proposition~7.8]{MR3272068}.
	
	\cref{75132773} Since $S(P)$ is a right $\Lambda_{\kp}$-module, $H_{\kp}$ is a nonzero $\Lambda_{\kp}$-submodule of $S(P)$. Thus $H_{\kp}=S(P)$ by \cref{13184618}.
	
	\cref{20391903} Let $0\neq L'\subset L\subset S(P)$ be $\Lambda$-submodules. Then by \cref{75132773}, $(S(P)/L)_{\kp}=0$ and $(L'/L)_{\kp}=0$. Thus, again by \cref{75132773}, $L'/L$ is not isomorphic to any nonzero $\Lambda$-submodule of $S(P)$. This means that $S(P)$ is a monoform right $\Lambda$-module.
\end{proof}

\begin{Convention}\label{74571009}
	For a noetherian $R$-algebra $\Lambda$ and $P\in\Spec\Lambda$, we always take $S(P)$ as the fixed representative of $\widetilde{P}\in\ASpec(\Mod\Lambda)$ for simplicity.
\end{Convention}

\begin{Proposition}\label{27923937}
	Let $\Lambda$ be a noetherian $R$-algebra. Let $P\in\Spec\Lambda$ and $\kp:=P\cap R$.
	\begin{enumerate}
		\item\label{95577757} For every integer $i\geq 0$, there is an isomorphism
		\begin{equation*}
			\Ext_{\Lambda}^{i}(\widetilde{P},-)=\varinjlim_{0\neq H\subset S(P)}\Ext_{\Lambda}^{i}(H,-)\isoto\Ext_{\Lambda_{\kp}}^{i}(S(P),(-)_{\kp})
		\end{equation*}
		of functors $\Mod\Lambda\to\Mod\bbZ$, induced from
		\begin{equation*}
			\Ext_{\Lambda}^{i}(H,-)\xrightarrow{(-)_{\kp}}\Ext_{\Lambda_{\kp}}^{i}(H_{\kp},(-)_{\kp})\isoto\Ext_{\Lambda_{\kp}}^{i}(S(P),(-)_{\kp})
		\end{equation*}
		for nonzero monoform $\Lambda$-submodules $H\subset S(P)$.
		\item\label{91194415} There is an isomorphism
		\begin{equation*}
			k(\widetilde{P})=\Hom_{\Lambda}(\widetilde{P},S(P))\isoto\End_{\Lambda_{\kp}}(S(P))
		\end{equation*}
		of skew fields, which is the isomorphism of \cref{95577757} applied to $S(P)$.
		
		We identify $k(\widetilde{P})$ with $\End_{\Lambda_{\kp}}(S(P))$ using this isomorphism. The isomorphism in \cref{95577757} can be regarded as that of functors $\Mod\Lambda\to\Mod k(\widetilde{P})$.
	\end{enumerate}
\end{Proposition}

\begin{proof}
	\cite[Proposition~7.8 and Theorem~7.10 (1)]{MR3272068}.
\end{proof}

\begin{Theorem}[\cref{09558590}]\label{89835504}
	Let $\Lambda$ be a noetherian $R$-algebra. Let $i\geq 0$ be an integer and $P,Q\in\Spec\Lambda$.
	\begin{enumerate}
		\item\label{34970095} If $P\cap R=Q\cap R=:\kp$, then $\varExt_{\Lambda}^{i}(\widetilde{P},\widetilde{Q})$ is eventually constant, and there is an isomorphism
		\begin{equation*}
			\Ext_{\Lambda}^{i}(\widetilde{P},\widetilde{Q})\isoto\Ext_{\Lambda_{\kp}}^{i}(S(P),S(Q))
		\end{equation*}
		of right $(k(\widetilde{Q})^{\op}\otimes_{R}k(\widetilde{P}))$-modules, induced from the isomorphism in \cref{27923937} \cref{95577757}.
		\item\label{58371959} If $P\cap R\neq Q\cap R$, then
		\begin{equation*}
			\Ext_{\Lambda}^{i}(\widetilde{P},\widetilde{Q})=0\quad\text{and}\quad\mathrm{D}_{\widetilde{P}}\varExt_{\Lambda}^{i}(\widetilde{P},\widetilde{Q})=0.
		\end{equation*}
		In particular, $\widetilde{Q}\notin\Epi_{i}(\widetilde{P})$ and $\widetilde{Q}\notin\Const_{i}(\widetilde{P})$.
	\end{enumerate}
\end{Theorem}

\begin{proof}
	\cref{34970095} Let $H\subset S(Q)$ be a nonzero $\Lambda$-submodule. Then by \cref{27923937}, we have the isomorphisms
	\begin{equation*}
		\Ext_{\Lambda}^{i}(\widetilde{P},H)\isoto\Ext_{\Lambda_{\kp}}^{i}(S(P),H_{\kp})\isoto\Ext_{\Lambda_{\kp}}^{i}(S(P),S(Q))
	\end{equation*}
	of right $k(\widetilde{P})$-modules. Thus we obtain the eventual constancy and the desired isomorphism of right $k(\widetilde{P})$-modules. It is straightforward to see that it is also an isomorphism of left $k(\widetilde{Q})$-modules.
	
	\cref{58371959} Set $\kp:=P\cap R$ and $\kq:=Q\cap R$. Let $H\subset S(Q)$ be a nonzero $\Lambda$-submodule.
	
	Assume $\kq\not\subset\kp$ and take $a\in\kq\setminus\kp$. Then $S(Q)Q=0$ and $a\in Q$ imply $H_{\kp}a=0$. Since $a$ is invertible in $\Lambda_{\kp}$, we have $H_{\kp}=0$. Thus $\Ext_{\Lambda}^{i}(\widetilde{P},H)=\Ext_{\Lambda_{\kp}}^{i}(S(P),H_{\kp})=0$ and the claims follow.
	
	Assume $\kp\not\subset\kq$ and take $b\in\kp\setminus\kq$. Since $b$ is invertible in $\Lambda_{\kq}$, it acts bijectively on $S(Q)$ and injectively on $H$. Thus the multiplication of $b$ is factorized as $H\isoto Hb\into H$ and induces a commutative diagram
	\begin{equation*}
		\begin{tikzcd}[column sep=small]
			\Ext_{\Lambda_{\kp}}^{i}(S(P),H_{\kp})\ar[dr,"\cong"']\ar[rr,"\cdot b"] & & \Ext_{\Lambda_{\kp}}^{i}(S(P),H_{\kp})\rlap{.} \\
			& \Ext_{\Lambda_{\kp}}^{i}(S(P),(Hb)_{\kp})\ar[ur,"f"'] &
		\end{tikzcd}
	\end{equation*}
	Since $b$ annihilates $S(P)$, the horizontal map is zero. So the map $f$, induced from the inclusion $Hb\into H$, is also zero. This concludes that, for every nonzero $\Lambda$-submodule $H\subset S(Q)$, there exists a nonzero $\Lambda$-submodule $H'\subset H$ such that the map
	\begin{equation*}
		\Ext_{\Lambda}^{i}(\widetilde{P},H')\to\Ext_{\Lambda}^{i}(\widetilde{P},H)
	\end{equation*}
	induced from the inclusion $H'\into H$ is zero. Since $\Ext_{\Lambda}^{i}(\widetilde{P},\widetilde{Q})$ (\resp $\mathrm{D}_{\widetilde{P}}\varExt_{\Lambda}^{i}(\widetilde{P},\widetilde{Q})$) is defined as an inverse limit (\resp a direct limit of the duals) of such maps, it should be zero.
\end{proof}

\begin{Corollary}\label{01099263}
	Let $R$ be a commutative noetherian ring. Let $i\geq 0$ be an integer and $\kp,\kq\in\Spec R$. Then
	\begin{equation*}
		\Ext_{R}^{i}(\widetilde{\kp},\widetilde{\kq})\cong
		\begin{cases}
			\Ext_{R_{\kp}}^{i}(k(\kp),k(\kp)) & \text{if $\kp=\kq$,} \\
			0 & \text{if $\kp\neq\kq$.}
		\end{cases}
	\end{equation*}
	If $\kp=\kq$, then $\varExt_{R}^{i}(\widetilde{\kp},\widetilde{\kp})$ is eventually constant.
\end{Corollary}

\begin{proof}
	Apply \cref{89835504} to $\Lambda=R$. For every $\kp\in\Spec R$, $S(\kp)$ is a simple $R_{\kp}$-module. Thus it is isomorphic to the residue field $k(\kp):=R_{\kp}/\kp R_{\kp}$.
\end{proof}

Consequently, we recover the following result of Gabriel:

\begin{Corollary}[{Gabriel \cite[Proposition~10 in p.~428]{MR0232821}}]\label{31584696}
	Let $R$ be a commutative noetherian ring.
	\begin{enumerate}
		\item\label{51857501} Every localizing subcategory of $\Mod R$ is closed under injective envelopes.
		\item\label{39903501} Every Serre subcategory of $\mod R$ is closed under essential extensions.
	\end{enumerate}
\end{Corollary}

\begin{proof}
	The claims follow from \cref{07240397,01099263}.
\end{proof}

\begin{Example}\label{62734803}
	Let $R$ be a commutative noetherian ring and consider the noetherian $R$-algebra
	\begin{equation*}
		\Lambda:=\begin{pmatrix}R & 0 \\ R & R\end{pmatrix}
	\end{equation*}
	as in \cite[Example~7.11]{MR3272068}, where $R$ is identified with the subring of all diagonal matrices with equal entries. All prime ideals of $\Lambda$ are of the form $P_{i}(\kp)$ with $i=1,2$ and $\kp\in\Spec R$, where
	\begin{equation*}
		P_{1}(\kp):=\begin{pmatrix}\kp & 0 \\ R & R\end{pmatrix}\quad\text{and}\quad P_{2}(\kp):=\begin{pmatrix}R & 0 \\ R & \kp\end{pmatrix}.
	\end{equation*}
	Note that $P_{i}(\kp)\cap R=\kp$. By \cref{63482658},
	\begin{equation*}
		\ASpec(\Mod\Lambda)=X_{1}\amalg X_{2},
	\end{equation*}
	where the right-hand side is the disjoint union of the two topological spaces
	\begin{equation*}
		X_{i}:=\setwithcondition{\widetilde{P_{i}(\kp)}}{\kp\in\Spec R},
	\end{equation*}
	and there are homeomorphisms $f_{i}\colon\ASpec(\Mod R)\isoto X_{i}$ given by $\overline{R/\kp}\mapsto\widetilde{P_{i}(\kp)}$. By \cref{09809844}, all localizing subcategories of $\Mod\Lambda$ are of the form
	\begin{equation}\label{99918584}
		\ASupp^{-1}(f_{1}(\Phi_{1})\cup f_{2}(\Phi_{2})),
	\end{equation}
	where $\Phi_{i}$ are open subsets of $\ASpec(\Mod R)$.
	
	It follows from \cref{89835504} that
	\begin{equation*}
		\Ext_{\Lambda}^{1}(\widetilde{P_{i}(\kp)},\widetilde{P_{j}(\kq)})=0\quad\text{if}\quad\kp\neq\kq
	\end{equation*}
	and
	\begin{equation*}
		\Ext_{\Lambda}^{1}(\widetilde{P_{i}(\kp)},\widetilde{P_{j}(\kp)})\cong\Ext_{\Lambda_{\kp}}^{1}(S(P_{i}(\kp)),S(P_{j}(\kp))).
	\end{equation*}
	We compute the last extension group. By \cite[Proposition~2.5.5]{MR1898632}, for every $\kp\in\Spec R$, the indecomposable injective right $\Lambda$-modules $I(P_{1}(\kp))$ and $I(P_{2}(\kp))$ are direct summands of $\Hom_{R}(\Lambda,E_{R}(R/\kp))$. It is easy to see that there is a decomposition
	\begin{equation*}
		\Hom_{R}(\Lambda,E_{R}(R/\kp))\cong\begin{pmatrix}E_{R}(R/\kp) & E_{R}(R/\kp)\end{pmatrix}\oplus\frac{\begin{pmatrix}E_{R}(R/\kp) & E_{R}(R/\kp)\end{pmatrix}}{\begin{pmatrix}E_{R}(R/\kp) & 0\end{pmatrix}}
	\end{equation*}
	into indecomposable injective right $\Lambda$-modules (which are unique up to isomorphism; see \cite[Proposition~2.7 (1)]{MR0099360}), where the action of $\Lambda$ in the right-hand side is defined by matrix multiplication. By \cref{63542092} \cref{13184618}, $S(P_{i}(\kp))$ is the unique simple $\Lambda_{\kp}$-submodule of $I(P_{i}(\kp))$ and it has to be annihilated by $P_{i}(\kp)$, so it follows that
	\begin{equation*}
		I(P_{1}(\kp))\cong\begin{pmatrix}E_{R}(R/\kp) & E_{R}(R/\kp)\end{pmatrix},\quad I(P_{2}(\kp))\cong\frac{\begin{pmatrix}E_{R}(R/\kp) & E_{R}(R/\kp)\end{pmatrix}}{\begin{pmatrix}E_{R}(R/\kp) & 0\end{pmatrix}},
	\end{equation*}
	and
	\begin{equation*}
		S(P_{1}(\kp))\cong\begin{pmatrix}k(\kp) & 0\end{pmatrix},\quad S(P_{2}(\kp))\cong\frac{\begin{pmatrix}k(\kp) & k(\kp)\end{pmatrix}}{\begin{pmatrix}k(\kp) & 0\end{pmatrix}}.
	\end{equation*}
	Since the inclusion $S(P_{j}(\kp))\into I(P_{j}(\kp))$ is an essential monomorphism in $\Mod\Lambda_{\kp}$ and $S(P_{i}(\kp))$ is a simple $\Lambda_{\kp}$-module, the inclusion induces an isomorphism
	\begin{equation*}
		\Hom_{\Lambda_{\kp}}(S(P_{i}(\kp)),S(P_{j}(\kp)))\isoto\Hom_{\Lambda_{\kp}}(S(P_{i}(\kp)),I(P_{j}(\kp))).
	\end{equation*}
	Hence the short exact sequence
	\begin{equation*}
		0\to S(P_{j}(\kp))\to I(P_{j}(\kp))\to I(P_{j}(\kp))/S(P_{j}(\kp))\to 0
	\end{equation*}
	induces an isomorphism
	\begin{equation*}
		\Hom_{\Lambda_{\kp}}(S(P_{i}(\kp)),I(P_{j}(\kp))/S(P_{j}(\kp)))\isoto\Ext_{\Lambda_{\kp}}^{1}(S(P_{i}(\kp)),S(P_{j}(\kp))).
	\end{equation*}
	By computing the left-hand side, we obtain
	\begin{equation*}
		\Ext_{\Lambda}^{1}(\widetilde{P_{i}(\kp)},\widetilde{P_{j}(\kp)})\cong
		\begin{cases}
			\Hom_{R_{\kp}}(k(\kp),E_{R}(R/\kp)/k(\kp))\cong\Ext_{R_{\kp}}^{1}(k(\kp),k(\kp)) & \text{if $i=j$,} \\
			\Hom_{R_{\kp}}(k(\kp),k(\kp))\cong k(\kp) & \text{if $i=2$ and $j=1$,} \\
			0 & \text{if $i=1$ and $j=2$.}
		\end{cases}
	\end{equation*}
	Therefore, by \cref{30980494}, the localizing subcategory of the form \cref{99918584} is closed under injective envelopes if and only if $\Phi_{1}\subset\Phi_{2}$.
\end{Example}



\end{document}